\input amstex
\documentstyle{amsppt}

\centerline{\bf ON THE NUMBER OF $\bold{8}$-CYCLES FOR TWO }

\smallskip

\centerline{\bf PARTICULAR REGULAR TOURNAMENTS OF ORDER $\bold{N}$}

\smallskip

\centerline{\bf WITH DIAMETRICALLY OPPOSITE  LOCAL PROPERTIES}

\medskip

\centerline{\bf S.V. Savchenko}

\medskip

\centerline{L.D. Landau Institute for Theoretical Physics, Russian Academy of Sciences}
\centerline{Kosygin str. 2, Moscow 119334, Russia
\footnote[]{E-mail address: savch$\@$itp.ac.ru}}

\medskip

\centerline{\sl Dedicated to Professor John Mackey on the occasion of his 60th birthday}

\bigskip

{\bf Abstract}

For a regular tournament $T$ of order $n,$
denote by $c_{8}(T)$ the number of cycles of length $8$ in $T.$
Let $DR_{n}$ be a doubly-regular tournament of order $n\equiv 3\mod4$
(so, the out-sets and in-sets of its vertices are also regular and hence,
contain the maximum possible number of cyclic triples)
and $RLT_{n}$ be the unique regular locally transitive tournament
of (odd) order $n$
(so, the out-sets and in-sets of its vertices are transitive and hence,
contain no cyclic triples, at all).
Some arguments based on the spectral properties of tournaments
allow us to suggest that $c_{8}(T) \le c_{8}(RLT_{n}),$
where $n$ is sufficiently large.
This restriction on $n$ is essential because our computer processing
of B. McKay's file of tournaments
implies that for $n=9,11,13,$
the maximum of $c_{8}(T)$ is attained at
tournaments with regular structure of the out and in-sets of their vertices.
In the present paper, we show that $c_{8}(DR_{n})$ does not depend on
a particular choice of $DR_{n}$ and determine expressions for
$c_{8}(DR_{n})$ and $c_{8}(RLT_{n}).$ They are both polynomials
of degree $8$ in $n.$ Comparing $c_{8}(DR_{n})$
with $c_{8}(RLT_{n})$ yields the inequality $c_{8}(DR_{n})>c_{8}(RLT_{n})$
for $11\le n\le 35,$ while $c_{8}(RLT_{n}) > c_{8}(DR_{n})$ for $n\ge 39.$
This allows us to treat the value $n=39$ as the point
of phase transition in the local properties of maximizers and
minimizers of $c_{8}(T)$ in the class of regular tournaments of order $n.$

\smallskip

Bibliography: 16 titles.

\smallskip


\smallskip


\bigskip

\noindent {\bf \S 1. Introduction}

\medskip

A {\sl tournament} $T$ of order $n$ is an orientation of the complete
graph $K_{n}.$ If a pair $(i,j)$ of vertices is an arc in $T,$  then
we say that the vertex $i$ dominates the vertex $j$ and write $i\to j.$
The {\sl in-set} $N^{-}(i)$ is the set of vertices dominating $i$ in $T.$
In turn, the {\sl out-set} $N^{+}(i)$ is the set of vertices dominated by
$i$ in $T.$ Denote by $\delta^{+}_{i}$ and $\delta^{-}_{i}$ their cardinalities. The
quantities $\delta^{+}_{i}$ and $\delta^{-}_{i}$ are  called the {\sl out-degree} and
{\sl in-degree} of the vertex $i,$ respectively. Clearly, $N^{+}(i)$ and
$N^{-}(i)$ induce subtournaments of orders $\delta^{+}_{i}$ and $\delta^{-}_{i},$ respectively.
In the sequel, we denote them by the same symbols $N^{+}(i)$ and $N^{-}(i).$

Let $c_{m}(T)$ be the number of cycles of length $m$ (or, merely, $m$-cycles)
in $T.$ The formula obtained for $c_{3}(T)$ by Kendall and Babington Smith in
[6] can be rewriten as follows:
$$2c_{3}(T)=\frac{n(n-1)(2n-1)}{6}-\sum\limits_{i=1}^{n}\bigl(\delta_{i}^{+}\bigr)^{2}=
\frac{n(n-1)(2n-1)}{6}-\sum\limits_{i=1}^{n}\bigl(\delta_{i}^{-}\bigr)^{2}.$$
In particular, if we take $\delta_{i}^{+}=i-1,$ we get $c_{3}(T)=0.$
For any given positive integer $n,$ there exists exactly one (up to
isomorphism) tournament of order $n$ having no $3$-cycles.
Indeed, if $1,...,n$ is a hamiltonian path (whose existence is quaranteed
by Redei's theorem) and $j\to i$ for $j>i+1,$ then $i,i+1,...,j,i$ is a
cycle. By the Moon theorem, the (strong) subtournament induced by the vertex-set
$\{i,...,j\}$ also contains a $3$-cycle. Hence, $i\to j$ if $j>i.$ Obviously,
this tournament also has the transitive property, i.e. the relation $i\to k\to j$
for three vertices $i,k,$ and $j$ implies $i\to j.$ By this reason, it is often
called the {\sl transitive} tournament of order $n.$ In the sequel, we denote
it by $TT_{n}(1,...,n)$ or, simply, by $TT_{n}.$

Note that $\sum\limits_{i=1}^{n}\delta_{i}^{+}=\frac{n(n-1)}{2}.$
Under this restriction, the minimum of the sum
$\sum\limits_{i=1}^{n}\bigl(\delta_{i}^{+}\bigr)^{2}$ is attained if and only if
$\delta_{i}^{+}=\frac{n-1}{2}$ for each $i.$ A tournament with this list of out-degrees
is called {\sl regular}. Obviously,  for this case, the order $n$ must be odd
and the in-degree of every vertex is also equal to $\frac{n-1}{2}.$ By this
reason, the quantity $\frac{n-1}{2}$ is called the {\sl semidegree} and is denoted by $\delta.$
A regular tournament of order $n$ exists for each odd $n\ge 1.$ So, for any
such $n,$ the maximum of $c_{3}(T)$ in the class ${\Cal T}_{n}$ of all tournaments
on $n$ vertices equals $\frac{n(n+1)(n-1)}{24}$ and
is achieved only for a regular tournament of order $n.$

A formula for $c_{4}(T)$ in the class ${\Cal T}_{n}$ is also well known.
For the case of a regular tournament $T$ of order $n,$ it can be presented
in the following form: $c_{4}(T)=$
$$\frac{n(n+1)(n-1)(n-3)}{48}-\sum\limits_{i=1}^{n}
c_{3}\bigl(N^{+}(i)\bigr)=
\frac{n(n+1)(n-1)(n-3)}{48}-\sum\limits_{i=1}^{n}
c_{3}\bigl(N^{-}(i)\bigr).$$

This allows us to get a two-sided bound on $c_{4}(T)$
in the class ${\Cal R}_{n}$
of regular tournaments of order $n.$ In particular, an upper bound can be
obtained if we assume that $c_{3}\bigl(N^{+}(i)\bigr)=0$ (and hence,
$c_{3}\bigl(N^{-}(i)\bigr)=0$) for each $i.$ This means that for each $i,$
both $N^{+}(i)$ and $N^{-}(i)$ are transitive.
It is not difficult to show that for each odd $n\ge 1,$
there exists exactly one tournament with this property (see [2]).
We call it the regular {\sl locally transitive} tournament of order $n$
and denote it by $RLT_{n}.$
It can be represented as the tournament
on the ring
${\Bbb Z}_{n}=\{0,...,n-1\}$ of residues modulo $n$
for which a pair $(i,j)$ is an arc if and only if
$j-i\in \{1,...,\frac{n-1}{2}\},$ where the subtraction is
taken modulo $n.$

In turn, a lower bound
in the class ${\Cal R}_{2\delta+1}$
of regular tournaments of semidegree $\delta$
(and hence, of order $2\delta+1$)
can be obtained if we assume that $c_{3}\bigl(N^{+}(i)\bigr)=
\frac{\delta(\delta+1)(\delta-1)}{24}$ (and hence,
$c_{3}\bigl(N^{-}(i)\bigr)=\frac{\delta(\delta+1)(\delta-1)}{24}$)
for each $i.$ This means that for each $i,$
both $N^{+}(i)$ and $N^{-}(i)$ are regular tournaments of order $\delta.$
A tournament of order $n$ with this property is called
{\sl doubly-regular} and is denoted by $DR_{n}.$

A well-known example of $DR_{p},$ where $p$ is a prime of the form $4r+3,$
is the tournament on (the field) ${\Bbb Z}_{p}$ with $(i,j)$ is an arc
if and only if $j-i$ is a non-zero square, called a quadratic residue.
(By this reason, such a $DR_{p}$ is called the {\sl quadratic residue
tournament} and is denoted by $QR_{p}$.)
In particular, for $n=3,7,$ and $11,$ a unique $DR_{n}$ is also $QR_{n},$
while for $n=15,$ there exist exactly two $DR_{n}$'s, which are not
vertex-transitive, at all. There are many
methods for construction of doubly-regular tournaments of other orders.
Nevertheless, the problem of the existence of $DR_{4t+3}$
for each $t\ge 0$ is open up to now, while according to a common opinion,
it exists for any possible order.

The results on $c_{4}(T)$ presented above (they were first obtained in
[1,16] and [8])
imply that the inequality
$$c_{4}(DR_{n})\le c_{4}(T)\le c_{4}(RLT_{n})$$
holds for any regular tournament $T$ of odd order $n.$
\footnote[1]{Note that while $DR_{n}$ can exist only for $n\equiv 3\mod 4,$
the expression for $c_{4}(DR_{n})$ is a polynomial of degree $4$ in $n$ and hence,
is defined for each integer value of $n$. In fact, according to [1] and [16],
for $n\equiv 1\mod 4,$ one can replace $c_{4}(DR_{n})$ in the left-hand side
by $c_{4}(NDR_{n}),$ where $NDR_{n}$ is a nearly-doubly-regular tournament
of order $n$ (by definition, the out-set (and in-set) of each of its vertices
is near regular, i.e. it is obtained from a regular tournament
by removing exactly one vertex).} (Note that the maximum number of $4$-cycles
in the class of all tournaments of odd order $n$ is also attained at $RLT_{n}$.)
However, since
$$c_{5}(T)+2c_{4}(T)=\frac{n(n-1)(n+1)(n-3)(n+3)}{160},$$
the converse inequality
$$c_{5}(RLT_{n})\le c_{5}(T)\le c_{5}(DR_{n})$$
holds for any regular tournament $T$ of order $n$ (see [12];
for the class ${\Cal T}_{n},$ where $n$ is odd, the inequality
$c_{5}(T)\le c_{5}(DR_{n})$ is proved in [15] with the use of the Komarov-Mackey
formula for $c_{5}(T)$ obtained in [7]).
In [12], we also show that $c_{6}(T)\le c_{6}(DR_{n})$
in the class
of all regular tournaments of order $n\ge 7.$
(The conjecture that $c_{6}(T)\le c_{6}(DR_{n})$ for each $T\in {\Cal T}_{n}$
is open up to now and in our opinion, it is a very difficult problem.)
Moreover, our computer search and
some pure spectral arguments allow us to
suggest that the minimum number
of $6$-cycles in the class ${\Cal R}_{n}$
is attained at $RLT_{n}$ for arbitrary odd $n\ge 7.$

In [13], we also determine expressions for
$c_{7}(DR_{n})$ and $c_{7}(RLT_{n}).$ They imply that
$c_{7}(DR_{n}) > c_{7}(RLT_{n})$ for each odd $n\ge 7.$
Moreover, some additional arguments allow us to suggest that
the two-sided bound $c_{7}(RLT_{n}) \le c_{7}(T)\le c_{7}(DR_{n})$
holds for each odd $n\ge 9.$
So, the results obtained imply that
$DR_{n}$ and $RLT_{n}$ are possible extreme points
of $c_{m}(T)$ in this class.
Hence,
it is important to determine $c_{m}(DR_{n})$ and $c_{m}(RLT_{n})$
at least for sufficiently small values of $m$
or compare them at least for sufficiently large values of $n.$
In [13] and [14],
we show that
$c_{m}(DR_{n})> c_{m}(RLT_{n})$ for $m\equiv 1,2,3 \mod 4$
and
$c_{m}(RLT_{n})> c_{m}(DR_{n})$ for $m\equiv 0 \mod 4$
if $n$ is sufficiently large.

The results obtained in the present paper imply that
the condition on the value of $n$
is essential at least in the case of $m=8.$
The paper is organized as follows. In Section 2,
an expression for $c_{8}(DR_{n})$ is determined on the basis of
the inclusion-exclusion principle already used in [12] and [13]
for determining $c_{6}(DR_{n})$ and $c_{7}(DR_{n}),$
respectively. As a by-product result, we also show that
any two arcs in a doubly-regular tournament
lie on the same number of $8$-cycles. Hence, the following
important conjecture first presented in [12] is confirmed for $m=8$:

\smallskip

{\bf Conjecture A}. {\sl For each $m\le 9,$
any two doubly-regular tournaments of order $n\equiv 3\mod 4$ have
the same number $c_{m}(DR_{n})$ of $m$-cycles.
Moreover, for these values of $m,$
any arc in a doubly-regular tournament
lies on $\frac{2mc_{m}(DR_{n})}{n(n-1)}$ cycles of length $m.$}

Our computer processing of
B. McKay's files of tournaments shows that
the restriction \lq\lq $m\le 9$\rq\rq
on the length $m$ is essential.
Note that for $m=3,$ $m=4,5,$ $m=6,$ and $m=7,$
Conjecture A was confirmed in [10], [9,11], [12], and [13], respectively.
The remaining case $m=9$ is more difficult and
will be considered in the separate paper
"On the number of $9$-cycles in a doubly-regular tournament",
which is now in preparation.
Moreover, it was also shown in [10], [1,16],  and [12]
that
\lq\lq any arc lies on $\frac{2mc_{m}(DR_{n})}{n(n-1)}$ cycles of length $m$\rq\rq
is a characterizing property of $DR_{n}$ in the class of tournaments
of order $n\ge m$ for $m=3,$ $m=4$ and $m=5,6,$ respectively. For $m=7,8,9,$
this problem is still open.

In Section 3, we obtain an expression for $c_{8}(RLT_{n})$
with the use of the method developed in [13]
for determining $c_{m}(RLT_{n})$
in the case of arbitrary $m.$ Comparing $c_{8}(DR_{n})$ with $c_{8}(RLT_{n})$
(they are polynomials of degree $8$ in $n$)
yields the strict inequality $c_{8}(RLT_{n})> c_{8}(DR_{n})$ for
all $n\ge 39,$ while $c_{8}(RLT_{n})< c_{8}(DR_{n})$ for $9\le n\le 35.$
This result, some spectral arguments presented in [13],
and the computer processing of B. McKay's file of tournaments
allow us to conjecture the existence of $N(8)\ge 39$
such that
$c_{8}(DR_{n})\le c_{8}(T)\le c_{8}(RLT_{n})$
for any regular tournament $T$ of order $n\ge N(8).$

Note that for determining $c_{8}(RLT_{n})$ with the use of
the method presented in [13], we need to know
the number $c_{7}(RLT_{n})$ of $7$-cycles in $RLT_{n}$
and the number $c_{8}^{(2)}(RLT_{n},0)$
of closed walks of length $8$ starting and ending at $0,$
containing exactly one copy of $0$ as an intermediate vertex,
and containing no other repeated vertices.
An expression for $c_{7}(RLT_{n})$ was first determined in [13].
We also reproved therein the known results on $c_{4}(RLT_{4}),$
$c_{5}(RLT_{n}),$ and $c_{6}(RLT_{n})$  obtained
in  [8], [3], [12], respectively.
The determining of $c_{8}^{(2)}(RLT_{n},0)$
is a much more difficult problem.
In Section 4, we solve it.
Though the most of
the corresponding calculations are routine, they are necessary
and we are not sure that one can significantly reduce them.

\bigskip

\noindent {\bf \S 2. The number of $\bold{8}$-cycles in $\bold{DR_{n}}$}

\medskip

Let $c_{m}(T,ji)$ be
the number of cycles of length $m$ including the arc $ji.$
Obviously, $c_{m}(T,ji)$  coincides with the number $p_{m-1}(T,i,j)$
of paths of length $m-1$ from
$i$ to $j.$
It is well known that
the number of all walks of length $m-1$
from $i$ to $j$ is equal
to the $ij$th entry $a_{ij}^{(m-1)}$
of the $(m-1)$th power of the adjacency matrix $A$ of $T.$
So, to determine $p_{m-1}(T,i,j),$
we need subtract from $a_{ij}^{(m-1)}$
the number $Corr_{m-1}(T,i,j)$
of walks of length $m-1$
from $i$ to $j$ in $T$
that are not paths.
In this section, we first determine the $ij$th entry
$a_{ij}^{(7)}$
of the $7$th power of the adjacency matrix $A$ of $DR_{4t+3}$
and then express
$Corr_{7}(DR_{4t+3},i,j)$ via the entries of the powers of the matrix $A.$

Note that the definition of $DR_{4t+3}$
implies that
any two vertices $i$ and $j$ of $DR_{4t+3}$
jointly dominate $t$ vertices, and
the number of vertices which
dominate both $i$ and $j$
is also equal to $t.$
Hence,
$a_{ij}^{(2)}=t$
and
$a_{ij}^{(2)}=t+1$
if $i\to j$ and $j\to i$ in $DR_{4t+3},$ respectively.
In particular, $a_{ii}^{(3)}=\sum\limits_{j\to i}a_{ij}^{(2)}=(2t+1)(t+1).$

It is also shown in [11] that the adjacency matrix $A$
of any $DR_{4t+3}$ satisfies the equation
$A^{3}=2tA^{2}+tA+(2t+1)(t+1)I,$
where $I$ is
the identity matrix, and, hence,
$$a^{(p)}_{ij}=2ta^{(p-1)}_{ij}+ta^{(p-2)}_{ij}+(2t+1)(t+1)a^{(p-3)}_{ij}$$
for each integer $p\ge 3.$
In particular, $a_{ij}^{(3)}=2ta_{ij}^{(2)}+ta_{ij}=
2t^{2}+t=t(2t+1)$ if $i\to j$  and
$a_{ij}^{(3)}=2ta_{ij}^{(2)}=
2t(t+1)$ if $j\to i$ in $DR_{4t+3}.$
This implies that $a_{ii}^{(4)}=(2t+1)\times 2t(t+1)$
for each vertex $i$ in $DR_{4t+3}.$
Moreover,
$$a_{ij}^{(4)}=2ta_{ij}^{(3)}+ta_{ij}^{(2)}=
2t\times 2t(t+1)+t\times (t+1)=
t(t+1)(4t+1)$$
if $ji$ is an arc in $DR_{4t+3}$
and hence, $a_{ii}^{(5)}=(2t+1)\times (4t+1)(t+1)t.$
Finally,

$$a_{ij}^{(5)}=2ta_{ij}^{(4)}+ta_{ij}^{(3)}+(2t+1)(t+1)a_{ij}^{(2)}=$$
$$2t^{2}(t+1)(4t+1)+2t^{2}(t+1)+(2t+1)(t+1)^{2}=$$
$$=(t+1)\Bigl(8t^{3}+2t^{2}+2t^{2}+(2t+1)(t+1)\Bigr)=(t+1)(2t+1)(4t^{2}+t+1),$$

$$a_{ij}^{(6)}=2ta_{ij}^{(5)}+ta_{ij}^{(4)}+(2t+1)(t+1)a_{ij}^{(3)}=
2t(t+1)(2t+1)(4t^{2}+t+1)+$$
$$t^{2}(t+1)(4t+1)+(2t+1)(t+1)2t(t+1)=t(t+1)\times $$
$$\Bigl((4t+2)(4t^{2}+t+1)
+t(4t+1)+2(t+1)(2t+1)\Bigr)=$$
$$t(t+1)\Bigl(2(2t+1)(4t^{2}+2t+2)+t(4t+1)\Bigr)
=t(t+1)\Bigl(16t^{3}+20t^{2}+13t+4\Bigr),$$

and hence,

$$a_{ij}^{(7)}=2ta_{ij}^{(6)}+ta_{ij}^{(5)}+(2t+1)(t+1)a_{ij}^{(4)}=
2t\times t(t+1)\Bigl(16t^{3}+20t^{2}+13t+4\Bigr)+$$
$$t\times (t+1)(2t+1)(4t^{2}+t+1)+(2t+1)(t+1)\times (4t+1)(t+1)t$$
$$=2t(t+1)\Bigl(16t^{4}+28t^{3}+23t^{2}+9t+1\Bigr).$$

Let us now determine $Corr_{7}(DR_{4t+3},i,j),$
where the pair $ji$ forms an arc in $DR_{4t+3}.$
Obviously, any walk of length $7$ that is not a path
has the form $\gamma=i,k,x,y,z,r,p,j$
in which at least one of the vertices $i,j,k,x,y$ is repeated.
Based on the inclusion-exclusion principle,
we will consequently consider all these five cases.

Suppose first that the vertex $i$ is
repeated in $\gamma.$ Then either $y=i,$ $z=i$ or $r=i.$
In the first case
$\gamma=i,k,x,i,z,r,p,j.$
The number of such walks from $i$ to $j$
is equal to the product of the number of closed walks (circuits)
of length $3$ starting and ending at $i$
and the number of walks (paths) of length $4$ from $i$ to $j.$
In turn, the number of the walks
of the form $\gamma=i,k,x,y,i,r,p,j$ equals the product of the number of
closed walks (cycles) of length $4$
starting and ending at $i$
and the number of walks (paths) of length $3$ from $i$ to $j.$
Finally, the number of the walks
of the form $\gamma=i,k,x,y,z,i,p,j$ equals the product of the number of
closed walks (cycles) of length $5$
starting and ending at $i$
and the number of walks (paths) of length $2$ from $i$ to $j.$
This can be reflected as follows:

\medskip

$i:$\ \ \ \ \ \ $\bold{i},*,*,\bold{i},*,*,*,\bold{j}$ \  \  \ \qquad \qquad $a_{ii}^{(3)}a_{ij}^{(4)}=
+(4t+1)(2t+1)(t+1)^{2}t$

\smallskip

\ \ \ \ \ \ \ \ \ $\bold{i},*,*,*,\bold{i},*,*,\bold{j}$ \  \  \   \qquad \qquad $a_{ii}^{(4)}a_{ij}^{(3)}=
+4(2t+1)(t+1)^{2}t^{2}$

\ \ \ \ \ \ \ \ \  $\bold{i},*,*,*,*,\bold{i},*,\bold{j}$ \  \  \   \qquad \qquad $a_{ii}^{(5)}a_{ij}^{(2)}=
+(4t+1)(2t+1)(t+1)^{2}t.$

\bigskip

In turn, if a walk $\gamma=i,k,x,y,z,r,p,j$ contains $j$ as
an intermediate vertex, then either $x=j$ or $y=j$ or $z=j.$ More precisely,

\medskip

$j:$ \ \ \ \ \ $\bold{i},*,\bold{j},*,*,*,*,\bold{j}$ \  \  \ \ \qquad \qquad $a_{ij}^{(2)}a_{jj}^{(5)}=
+(4t+1)(2t+1)(t+1)^{2}t$

\smallskip

\ \ \ \ \ \ \ \ \ $\bold{i},*,*,\bold{j},*,*,*,\bold{j}$ \  \  \  \qquad \qquad $a_{ij}^{(3)}a_{jj}^{(4)}=
+4(2t+1)(t+1)^{2}t^{2}$

\smallskip

\ \ \ \ \ \ \ \ \ $\bold{i},*,*,*,\bold{j},*,*,\bold{j}$ \  \  \ \qquad \qquad $a_{ij}^{(4)}a_{jj}^{(3)}=
+(4t+1)(2t+1)(t+1)^{2}t.$

\bigskip

But some of these walks
also repeat $i$
and hence, they must be subtracted from the previous sum:

\medskip

$\overset--\to{ji}:$\ \ \ \ $\bold{i},*,\bold{j},\bold{i},*,*,*,\bold{j}$\  \  \ \qquad  \qquad
$-a_{ij}^{(2)}a_{ij}^{(4)}=-(4t+1)(t+1)^{2}t$

\smallskip

\ \ \ \ \ \ \ \ \ $\bold{i},*,\bold{j},*,\bold{i},*,*,\bold{j}$\  \  \ \ \qquad \qquad
$-a_{ij}^{(2)}a_{ji}^{(2)}a_{ij}^{(3)}=-2(t+1)^{2}t^{2}$

\smallskip

\ \ \ \ \ \ \ \ \ $\bold{i},*,\bold{j},*,*,\bold{i},*,\bold{j}$\  \  \  \ \qquad \qquad
$-a_{ij}^{(2)}a_{ji}^{(3)}a_{ij}^{(2)}=-(2t+1)(t+1)^{2}t$

\smallskip

\ \ \ \ \ \ \ \ \ $\bold{i},*,*,\bold{j},\bold{i},*,*,\bold{j}$\  \  \  \ \qquad \qquad
$-a_{ij}^{(3)}a_{ij}^{(3)}=-4(t+1)^{2}t^{2}$

\smallskip

\ \ \ \ \ \ \ \ \ $\bold{i},*,*,\bold{j},*,\bold{i},*,\bold{j}$\  \  \  \ \qquad \qquad
$-a_{ij}^{(3)}a_{ji}^{(2)}a_{ij}^{(2)}=-2(t+1)^{2}t^{2}$

\smallskip

\ \ \ \ \ \ \ \ \ $\bold{i},*,*,*,\bold{j},\bold{i},*,\bold{j}$\  \  \  \ \qquad \qquad
$-a_{ij}^{(4)}a_{ij}^{(2)}=-(4t+1)(t+1)^{2}t.$

\medskip

Assume now that the vertex $k$ is repeated in a walk $\gamma=i,k,x,y,z,r,p,j.$
Obviously, either $z=k$ or  $r=k$ or $p=k.$
In the third case, we necessarily have $k\to j.$
For the first two cases, both orientations of the edge between $k$ and $j$
are possible. Hence,

\medskip

$k:$\ \ \ \ \ \ \ $\bold{i},k,*,*,k,*,*,\bold{j}$  \  \  \  \qquad \qquad \quad
$a_{kk}^{(3)}a_{\overset\curvearrowleft\to{ij}}^{(4)}=+(4t+1)(2t+1)(t+1)^{2}t$

\smallskip

\ \ \ \ \ \ \ \ \ \ $\bold{i},k,*,*,*,k,*,\bold{j}$  \  \  \ \qquad \qquad \quad
$a_{kk}^{(4)}a_{\overset\curvearrowleft\to{ij}}^{(3)}=+4(2t+1)(t+1)^{2}t^{2}$

\smallskip

\ \ \ \ \ \ \ \ \ \ $\bold{i},k,*,*,*,*,k,\bold{j}$  \  \  \  \qquad \qquad \quad
$a_{kk}^{(5)}a_{\overset\curvearrowleft\to{ij}}^{(2)}=+(4t+1)(2t+1)(t+1)^{2}t.$

\bigskip

Note that some of the walks repeating $k$
also contain $i$ and $j$ as a intermediate vertices.
Hence, we must subtract them from the sum obtained above:

\medskip

$\overset--\to{ki}:$\ \ \ \ $\bold{i},k,*,\bold{i},k,*,*,\bold{j}$\  \  \  \ \qquad  \qquad
$-a_{\overset\curvearrowleft\to{ki}}^{(2)}a_{\overset\curvearrowleft\to{ij}}^{(4)}=-(4t+1)(t+1)^{2}t$

\smallskip

\ \ \ \ \ \ \ \ \ $\bold{i},k,*,\bold{i},*,k,*,\bold{j}$\  \  \  \ \qquad\qquad
$-a_{\overset\curvearrowleft\to{ki}}^{(2)}
a_{\overset\curvearrowright\to{ik}}^{(2)}
a_{\overset\curvearrowleft\to{ij}}^{(3)}=-2(t+1)^{2}t^{2}$

\smallskip

\ \ \ \ \ \ \ \ \ $\bold{i},k,*,*,\bold{i},k,*,\bold{j}$\  \  \ \ \qquad \qquad
$-a_{\overset\curvearrowleft\to{ki}}^{(3)}a_{\overset\curvearrowleft\to{ij}}^{(3)}=-4(t+1)^{2}t^{2}$

\smallskip

\ \ \ \ \ \ \ \ \ $\bold{i},k,*,\bold{i},*,*,k,\bold{j}$\  \  \  \ \qquad \qquad
$-a_{\overset\curvearrowleft\to{ki}}^{(2)}a_{\overset\curvearrowright\to{ik}}^{(3)}
a_{\overset\curvearrowleft\to{ij}}^{(2)}=-(2t+1)(t+1)^{2}t$

\smallskip

\ \ \ \ \ \ \ \ \ $\bold{i},k,*,*,\bold{i},*,k,\bold{j}$\  \  \  \ \qquad \qquad
$-a_{\overset\curvearrowleft\to{ki}}^{(3)}
a_{\overset\curvearrowright\to{ik}}^{(2)}
a_{\overset\curvearrowleft\to{ij}}^{(2)}=
-2(t+1)^{2}t^{2}$

\smallskip

\ \ \ \ \ \ \ \ \ $\bold{i},k,*,*,*,\widehat{\bold{i},k,\bold{j}}$\  \  \  \ \qquad \qquad
$-a_{\overset\curvearrowleft\to{ki}}^{(4)}
a_{\overset\curvearrowleft\to{ij}}^{(2)}=
-(4t+1)(t+1)^{2}t$

\medskip

$\overset--\to{kj}:$ \ \ \ $\widehat{\bold{i},k,\bold{j}},*,k,*,*,\bold{j}$\  \  \  \qquad \qquad
$-a_{\overset\curvearrowleft\to{ij}}^{(2)}a_{\overset\curvearrowleft\to{jk}}^{(2)}
a_{\overset\curvearrowright\to{kj}}^{(3)}=-(2t+1)(t+1)^{2}t$

\smallskip

\ \ \ \ \ \ \ \ \  $\bold{i},k,*,\bold{j},k,*,*,\bold{j}$\  \  \  \qquad \qquad
$-\sum\limits_{i\to k\gets j} a_{kj}^{(2)}a_{kj}^{(3)}=-2(t+1)^{2}t^{2}$

\smallskip

\ \ \ \ \ \ \ \ \ $\widehat{\bold{i},k,\bold{j}},*,*,k,*,\bold{j}$\  \  \  \qquad \qquad
$-a_{\overset\curvearrowleft\to{ij}}^{(2)}
a_{\overset\curvearrowleft\to{jk}}^{(3)}
a_{\overset\curvearrowright\to{kj}}^{(2)}=-2(t+1)^{2}t^{2}$

\smallskip

\ \ \ \ \ \ \ \ \  $\bold{i},k,*,\bold{j},*,k,*,\bold{j}$\  \  \  \qquad \qquad
$-a_{kj}^{(2)}a_{jk}^{(2)}a_{\overset\curvearrowleft\to{ij}}^{(3)}
=-2(t+1)^{2}t^{2}$

\smallskip

\ \ \ \ \ \ \ \ \  $\bold{i},k,*,*,\bold{j},k,*,\bold{j}$\  \  \  \qquad \qquad
$-\sum\limits_{i\to k\gets j} a_{kj}^{(3)}a_{kj}^{(2)}=-2(t+1)^{2}t^{2}$

\smallskip

\ \ \ \ \ \ \ \ \  $\widehat{\bold{i},k,\bold{j}},*,*,*,k,\bold{j}$\  \  \  \qquad\qquad
$-a_{\overset\curvearrowleft\to{ij}}^{(2)}
a_{\overset\curvearrowleft\to{jk}}^{(4)}=-(4t+1)(t+1)^{2}t$

\smallskip

\ \ \ \ \ \ \ \ \  $\bold{i},k,*,\bold{j},*,*,k,\bold{j}$\  \  \  \qquad \qquad
$-a_{\overset\curvearrowright\to{kj}}^{(2)}
a_{\overset\curvearrowleft\to{jk}}^{(3)}
a_{\overset\curvearrowleft\to{ij}}^{(2)}
=-2(t+1)^{2}t^{2}$

\smallskip

\ \ \ \ \ \ \ \ \  $\bold{i},k,*,*,\bold{j},*,k,\bold{j}$\  \  \  \qquad \qquad
$-a_{\overset\curvearrowright\to{kj}}^{(3)}
a_{\overset\curvearrowleft\to{jk}}^{(2)}
a_{\overset\curvearrowleft\to{ij}}^{(2)}
=-(2t+1)(t+1)^{2}t.$

\bigskip

Obviously, the walks repeating both $i$ and $k$
are counting twice. Thus, we must add them to the sum obtained at the previous stage:

\medskip

$\overset+++\to{kij}:$ \ \ $\widehat{\bold{i},k,\bold{j}},\bold{i},k,*,*,\bold{j}$\  \  \  \qquad \qquad
$+a_{\overset\curvearrowleft\to{ij}}^{(2)}
a_{\overset\curvearrowright\to{kj}}^{(3)}=+(2t+1)(t+1)t$

\smallskip

\ \ \ \ \ \ \ \ \ \ $\widehat{\bold{i},k,\bold{j}},\bold{i},*,k,*,\bold{j}$\  \  \  \qquad \qquad
$+a_{\overset\curvearrowleft\to{ij}}^{(2)}
a_{\overset\curvearrowright\to{ik}}^{(2)}
a_{\overset\curvearrowright\to{kj}}^{(2)}=+(t+1)t^{2}$

\smallskip

\ \ \ \ \ \ \ \ \ \ $\widehat{\bold{i},k,\bold{j}},*,\bold{i},k,*,\bold{j}$\  \  \  \qquad \qquad
$+a_{\overset\curvearrowleft\to{ij}}^{(2)}
a_{\overset\curvearrowright\to{ji}}^{(2)}
a_{\overset\curvearrowright\to{kj}}^{(2)}=+(t+1)t^{2}$

\smallskip

\ \ \ \ \ \ \ \ \ \  $\bold{i},k,*,\bold{j},\bold{i},k,*,\bold{j}$\  \  \  \qquad \qquad
$+\sum\limits_{i\to k\to j}a_{kj}^{(2)}
a_{kj}^{(2)}+\sum\limits_{i\to k\gets j}a_{kj}^{(2)}
a_{kj}^{(2)}$

\qquad\qquad\qquad\qquad \qquad \qquad \qquad\quad $=+(2t+1)(t+1)t$

\smallskip

\ \ \ \ \ \ \ \ \ \ $\widehat{\bold{i},k,\bold{j}},\bold{i},*,*,k,\bold{j}$\  \  \  \qquad \qquad
$+a_{\overset\curvearrowleft\to{ij}}^{(2)}
a_{\overset\curvearrowright\to{ik}}^{(3)}=+(2t+1)(t+1)t$

\smallskip

\ \ \ \ \ \ \ \ \ \ $\widehat{\bold{i},k,\bold{j}},*,\bold{i},*,k,\bold{j}$\  \  \  \qquad \qquad
$+a_{\overset\curvearrowleft\to{ij}}^{(2)}
a_{\overset\curvearrowright\to{ji}}^{(2)}
a_{\overset\curvearrowright\to{ik}}^{(2)}=+(t+1)t^{2}$

\smallskip

\ \ \ \ \ \ \ \ \ \ $\bold{i},k,*,\bold{j},\bold{i},*,k,\bold{j}$\  \  \  \qquad \qquad
$+
a_{\overset\curvearrowleft\to{ij}}^{(2)}
a_{\overset\curvearrowright\to{kj}}^{(2)}
a_{\overset\curvearrowright\to{ik}}^{(2)}=+(t+1)t^{2}$

\smallskip

\ \ \ \ \ \ \ \ \ \ $\widehat{\bold{i},k,\bold{j}},*,*,\bold{i},k,\bold{j}$\  \  \  \qquad \qquad
$+
a_{\overset\curvearrowleft\to{ij}}^{(2)}
a_{\overset\curvearrowright\to{ji}}^{(3)}=+(2t+1)(t+1)t$

\smallskip

\ \ \ \ \ \ \ \ \ \ $\bold{i},k,*,\bold{j},*,\widehat{\bold{i},k,\bold{j}}$\  \  \ \qquad \qquad
$+
a_{\overset\curvearrowright\to{kj}}^{(2)}
a_{\overset\curvearrowright\to{ji}}^{(2)}
a_{\overset\curvearrowleft\to{ij}}^{(2)}
=+(t+1)t^{2}$

\smallskip

\ \ \ \ \ \ \ \ \ \ $\bold{i},k,*,*,\bold{j},\widehat{\bold{i},k,\bold{j}}$\  \  \  \qquad \qquad
$+a_{\overset\curvearrowright\to{kj}}^{(3)}
a_{\overset\curvearrowleft\to{ij}}^{(2)}
=+(2t+1)(t+1)t.$

\bigskip

Thus, summing all the expressions obtained above yields the following

\smallskip

$I=6*(4t+1)(2t+1)(t+1)^{2}t+3*4(2t+1)(t+1)^{2}t^{2}
-5*(4t+1)(t+1)^{2}t+$

$-2*4(t+1)^{2}t^{2}-4*(2t+1)(t+1)^{2}t-9*2(t+1)^{2}t^{2}+5*(2t+1)(t+1)t+$

$+5*(t+1)t^{2}=6(6t+1)(2t+1)(t+1)^{2}t-9(6t+1)(t+1)^{2}t+5(3t+1)(t+1)t$

$=3(4t-1)(6t+1)(t+1)^{2}t+5(3t+1)(t+1)t=2t(t+1)(36t^{3}+33t^{2}+3t+1).$

\smallskip

Consider now the $7$-walks $\gamma=i,k,x,y,z,r,p,j$
in which the third vertex $x$ is repeated and does not coincide with $j$:

\medskip

$x:$\ \ \ \ \ \ \ $\bold{i},*,x,*,*,x,*,\bold{j}$\  \  \  \qquad \qquad
$a_{xx}^{(3)}a_{\overset\curvearrowleft\to{ij}}^{(4)}=(4t+1)(2t+1)(t+1)^{2}t$

\smallskip

\ \ \ \ \ \ \ \ \ \ $\bold{i},*,x,*,*,*,x,\bold{j}$\  \  \  \qquad \qquad
$a_{xx}^{(4)}a_{\overset\curvearrowleft\to{ij}}^{(3)}=4(2t+1)(t+1)^{2}t^{2}$

\medskip

$\overset--\to{xi}:$\ \ \ \ \ $\bold{i},*,x,\bold{i},*,x,*,\bold{j}$\  \  \ \qquad \qquad
$-\sum\limits_{i\gets x\to j}a_{ix}^{(2)}a_{ix}^{(2)}a_{xj}^{(2)}-
\sum\limits_{i\gets x\gets j}a_{ix}^{(2)}a_{ix}^{(2)}a_{xj}^{(2)}$

\qquad\qquad\qquad\qquad\qquad\qquad \qquad \quad\quad $=-(2t+1)(t+1)^{2}t$

\smallskip

\ \ \ \ \ \ \ \ \ \ \ $\bold{i},*,x,*,\bold{i},x,*,\bold{j}$\  \  \  \qquad \qquad
$-a_{ix}^{(2)}a_{xi}^{(2)}a_{\overset\curvearrowleft\to{ij}}^{(3)}
=-2(t+1)^{2}t^{2}$

\smallskip

\ \ \ \ \ \ \ \ \ \ \ $\bold{i},*,x,\bold{i},*,*,x,\bold{j}$\  \  \  \qquad\qquad
$-\sum\limits_{i\gets x\to j}a_{ix}^{(2)}a_{ix}^{(3)}=-2(t+1)^{2}t^{2}$

\smallskip

\qquad  \ \ \ \ \ $\bold{i},*,x,*,\bold{i},*,x,\bold{j}$\  \  \ \qquad\qquad
$-a_{ix}^{(2)}a_{xi}^{(2)}
a_{\overset\curvearrowleft\to{ij}}^{(3)}=-2(t+1)^{2}t^{2}$

\smallskip

\qquad \ \ \ \ \ $\bold{i},*,x,*,*,\widehat{\bold{i},x,\bold{j}}$\  \  \  \qquad \qquad
$-a_{\overset\curvearrowright\to{ix}}^{(2)}
a_{\overset\curvearrowleft\to{xi}}^{(3)}
a_{\overset\curvearrowleft\to{ij}}^{(2)}
=-2(t+1)^{2}t^{2}$

\medskip

$\overset--\to{xj}:$ \ \ \ \ $\bold{i},*,x,\bold{j},*,x,*,\bold{j}$ \  \  \ \qquad \qquad
$-a_{\overset\curvearrowleft\to{ij}}^{(3)}a_{jx}^{(2)}a_{xj}^{(2)}
=-2(t+1)^{2}t^{2}$

\smallskip

\qquad \ \ \ \ $\bold{i},*,x,*,\bold{j},x,*,\bold{j}$ \  \  \ \qquad \qquad
$-\sum\limits_{i\to x\gets j}a_{ix}^{(2)}a_{xj}^{(2)}a_{xj}^{(2)}-
\sum\limits_{i\gets x\gets j}a_{ix}^{(2)}a_{xj}^{(2)}a_{xj}^{(2)}$

\qquad\qquad\qquad\qquad\qquad\qquad\qquad \qquad $=-(2t+1)(t+1)^{2}t$

\smallskip
\ \ \ \ \ \ \ \ \ \ $\bold{i},*,x,\bold{j},*,*,x,\bold{j}$\  \  \  \qquad \qquad
$-a_{\overset\curvearrowleft\to{ij}}^{(3)}
a_{\overset\curvearrowleft\to{jx}}^{(3)}=-4t^{2}(t+1)^{2}$

\smallskip

\qquad  \ \ \ \ $\bold{i},*,x,*,\bold{j},*,x,\bold{j}$\  \  \ \qquad \qquad
$-a_{xj}^{(2)}a_{jx}^{(2)}a_{\overset\curvearrowleft\to{ij}}^{(3)}
=-2t^{2}(t+1)^{2}$

\medskip

$\overset+++\to{xji}:$ \ \ $\bold{i},*,x,\bold{j},\bold{i},x,*,\bold{j}$\  \  \ \ \qquad \qquad
$+\sum\limits_{i\to x\to j}a_{ix}^{(2)}a_{xj}^{(2)}=+(t+1)t^{2}$

\smallskip

\ \ \ \ \ \ \ \ \ \ \ $\bold{i},*,x,\bold{j},\bold{i},*,x,\bold{j}$\  \  \  \qquad \qquad
$+\sum\limits_{i\to x\to j}
a_{ix}^{(2)}a_{ix}^{(2)}
+\sum\limits_{i\gets x\to j}
a_{ix}^{(2)}a_{ix}^{(2)}$

\qquad\qquad\qquad\qquad\qquad\qquad\qquad \qquad \quad $=+(2t+1)(t+1)t$

\smallskip

\qquad \ \ \ \ \ $\bold{i},*,x,\bold{j},*,\widehat{\bold{i},x,\bold{j}}$\  \  \ \qquad \qquad
$+a_{\overset\curvearrowright\to{ix}}^{(2)}
a_{\overset\curvearrowright\to{ji}}^{(2)}
a_{\overset\curvearrowleft\to{ij}}^{(2)}
=+(t+1)t^{2}$

\smallskip

\qquad \ \ \ \ \ $\bold{i},*,x,*,\bold{j},\widehat{\bold{i},x,\bold{j}}$\  \  \  \qquad \qquad
$+a_{\overset\curvearrowright\to{ix}}^{(2)}
a_{\overset\curvearrowright\to{xj}}^{(2)}
a_{\overset\curvearrowleft\to{ij}}^{(2)}
=+(t+1)t^{2}$

\medskip

$\overset--\to{xk}:$\ \ \ \ \ \ $\bold{i},k,x,*,k,x,*,\bold{j}$\  \  \ \qquad \qquad
$-a_{\overset\curvearrowleft\to{xk}}^{(2)}
a_{\overset\curvearrowleft\to{ij}}^{(4)}=-(4t+1)(t+1)^{2}t$

\ \ \ \ \ \ \ \ \ \ \ $\bold{i},k,x,*,k,*,x,\bold{j}$\  \  \ \qquad \qquad
$-a_{xk}^{(2)}a_{kx}^{(2)}
a_{\overset\curvearrowleft\to{ij}}^{(3)}=-2t^{2}(t+1)^{2}$

\smallskip

\qquad \ \ \ \ \ $\bold{i},k,x,*,*,k,x,\bold{j}$\  \  \ \qquad \qquad
$-a_{\overset\curvearrowleft\to{xk}}^{(3)}
a_{\overset\curvearrowleft\to{ij}}^{(3)}=-4t^{2}(t+1)^{2}$

\medskip

$\overset+++\to{xki}:$\ \ \   $\widehat{\bold{i},k,x},\bold{i},k,x,*,\bold{j}$\ \ \ \  \qquad \qquad
$+\sum\limits_{i\gets x\to j}a_{ix}^{(2)}
a_{xj}^{(2)}
+\sum\limits_{i\gets x\gets j}a_{ix}^{(2)}a_{xj}^{(2)}$

\qquad\qquad\qquad\qquad\qquad\qquad\qquad \qquad \quad $=+(2t+1)(t+1)t$

\ \ \ \ \ \ \ \ \ \ \ $\widehat{\bold{i},k,x},\bold{i},k,*,x,\bold{j}$\  \  \ \qquad \qquad
$+\sum\limits_{i\gets x\to j}
a_{\overset\curvearrowleft\to{ix}}^{(2)}
a_{\overset\curvearrowright\to{kx}}^{(2)}=+(t+1)t^{2}$

\smallskip

\quad \qquad   \ \ $\widehat{\bold{i},k,x},\bold{i},*,k,x,\bold{j}$\  \  \  \qquad\qquad
$+\sum\limits_{i\gets x\to j}
a_{\overset\curvearrowleft\to{ix}}^{(2)}
a_{\overset\curvearrowright\to{ik}}^{(2)}=+(t+1)t^{2}$

\smallskip

\quad\qquad   \ \  $\widehat{\bold{i},k,x},*,\bold{i},k,x,\bold{j}$\  \  \  \qquad\qquad
$+\sum\limits_{i\to x\to j}
a_{ix}^{(2)}a_{xi}^{(2)}+\sum\limits_{i\gets x\to j}a_{ix}^{(2)}
a_{xi}^{(2)}$

\qquad\qquad\qquad\qquad\qquad\qquad\qquad \qquad \quad $=+(2t+1)(t+1)t$

\medskip

$\overset+++\to{xkj}:$ \ \ $\bold{i},\widehat{k,x,\bold{j}},k,x,*,\bold{j}$\  \  \  \  \qquad \qquad
$+\sum\limits_{i\to k\gets j}
a_{\overset\curvearrowleft\to{kj}}^{(2)}
a_{\overset\curvearrowright\to{xj}}^{(2)}=+(t+1)t^{2}$

\smallskip

\ \ \ \ \ \ \ \ \ \ \ $\bold{i},\widehat{k,x,\bold{j}},k,*,x,\bold{j}$\  \  \ \qquad \qquad
$+\sum\limits_{i\to k\gets j}
a_{\overset\curvearrowleft\to{kj}}^{(2)}
a_{\overset\curvearrowright\to{kx}}^{(2)}=+(t+1)t^{2}$

\smallskip

\qquad \ \ \ \ \ $\bold{i},\widehat{k,x,\bold{j}},*,k,x,\bold{j}$\  \  \ \qquad \qquad
$+\sum\limits_{i\to k}
a_{kj}^{(2)}a_{jk}^{(2)}=+(2t+1)(t+1)t$

\smallskip

\qquad \ \ \ \ \ $\bold{i},k,x,*,\bold{j},\widehat{k,x,\bold{j}}$\  \  \  \qquad  \qquad
$+\sum\limits_{i\to k\gets j}
a_{\overset\curvearrowright\to{xj}}^{(2)}
a_{\overset\curvearrowleft\to{kj}}^{(2)}
=+(t+1)t^{2}$

\medskip

$\overset----\to{xkij}:$ \ $\bold{i},k,x,\bold{j},\bold{i},k,x,\bold{j}$
\ \ \qquad\qquad  $-a_{\overset\curvearrowleft\to{ij}}^{(3)}=-2t(t+1).$

\bigskip

Thus, the sum associated with the case where $r=x$ or $p=x$ has the form:

\smallskip

$II=(4t+1)(2t+1)(t+1)^{2}t+4(2t+1)(t+1)^{2}t^{2}-2(2t+1)(t+1)^{2}t$

$-22(t+1)^{2}t^{2}-(4t+1)(t+1)^{2}t+8(t+1)t^{2}+4(2t+1)(t+1)t-2t(t+1)$

$=(8t+1)(2t+1)(t+1)^{2}t-(30t+3)(t+1)^{2}t+(16t+2)(t+1)t=$

$(16t^{2}-20t-2)(t+1)^{2}+(16t+2)(t+1)t=(16t^{3}-4t^{2}-22t-2)(t+1)t$

$+(16t+2)(t+1)t=2t(8t^{2}-2t-3)(t+1)t=2t^{2}(t+1)(4t-3)(2t+1).$

\medskip

Finally, consider the case
when the fourth vertex $y$ is repeated.
Obviously,  this is possible iff
$p=y$ and hence, $y\to j.$ More precisely,

\medskip

$y:$\ \ \ \ \ \ \ $\bold{i},*,*,y,*,*,y,\bold{j}$\  \  \  \qquad \qquad
$a_{yy}^{(3)}a_{\overset\curvearrowleft\to{ij}}^{(4)}=+(4t+1)(2t+1)(t+1)^{2}t$

\medskip

$\overset--\to{yi}:$\ \ \ \ \ $\bold{i},*,*,y,\bold{i},*,y,\bold{j}$\  \  \ \qquad \qquad
$-\sum\limits_{i\gets y\to j}
a_{iy}^{(3)}a_{iy}^{(2)}=-2(t+1)^{2}t^{2}$

\smallskip

\qquad \ \ \ \ \ $\bold{i},*,*,y,*,\widehat{\bold{i},y,\bold{j}}$\  \  \  \qquad  \qquad
$-
a_{\overset\curvearrowright\to{iy}}^{(3)}
a_{\overset\curvearrowleft\to{yi}}^{(2)}
a_{\overset\curvearrowleft\to{ij}}^{(2)}
=-(2t+1)(t+1)^{2}t$

\medskip

$\overset--\to{yj}:$\ \ \ \ \ \ \ $\bold{i},*,*,y,\bold{j},*,y,\bold{j}$\  \  \  \qquad \qquad
$-a_{\overset\curvearrowleft\to{ij}}^{(4)}
a_{\overset\curvearrowleft\to{jy}}^{(2)}=-(4t+1)(t+1)^{2}t$

\medskip

$\overset+++\to{yji}:$\ \ \ \ \ $\bold{i},*,*,y,\bold{j},\widehat{\bold{i},y,\bold{j}}$\  \  \ \qquad \qquad
$+a_{\overset\curvearrowright\to{iy}}^{(3)}
a_{\overset\curvearrowleft\to{ij}}^{(2)}
=+(2t+1)(t+1)t$

\medskip

$\overset--\to{yk}:$\ \ \ \ \ \ $\bold{i},k,*,y,k,*,y,\bold{j}$\  \  \  \qquad \qquad
$-\sum\limits_{y\to j} \sum\limits_{y\to k\gets i}
a_{ky}^{(2)}a_{ky}^{(2)}=-(2t+1)(t+1)^{2}t$

\smallskip

\qquad \ \ \ \ \ \ \ $\bold{i},k,*,y,*,k,y,\bold{j}$\  \  \ \qquad  \qquad
$-a_{ky}^{(2)}a_{yk}^{(2)}
a_{\overset\curvearrowleft\to{ij}}^{(3)}=-2(t+1)^{2}t^{2}$

\medskip

$\overset+++\to{yki}:$\ \ \ \ \ $\bold{i},k,*,y,\widehat{\bold{i},k,y},\bold{j}$\  \  \ \qquad \qquad
$+\sum\limits_{i\gets y\to j}
a_{\overset\curvearrowright\to{ky}}^{(2)}
a_{\overset\curvearrowleft\to{iy}}^{(2)}
=+(t+1)t^{2}$

\medskip

$\overset+++\to{ykj}:$\ \ \ \ \ $\bold{i},k,*,y,\bold{j},\widehat{k,y,\bold{j}}$\  \  \  \qquad \qquad
$+\sum\limits_{i\to k\gets j}
a_{\overset\curvearrowright\to{ky}}^{(2)}
a_{\overset\curvearrowleft\to{kj}}^{(2)}
=+(t+1)t^{2}$

\medskip

$\overset--\to{yx}:$\ \ \ \ \ \ \ $\bold{i},*,x,y,*,x,y,\bold{j}$\  \  \ \qquad\qquad
$-a_{\overset\curvearrowleft\to{yx}}^{(2)}
a_{\overset\curvearrowleft\to{ij}}^{(4)}=-(4t+1)(t+1)^{2}t$

\medskip

$\overset+++\to{yxi}:$\ \ \ \ \ $\bold{i},*,x,y,\widehat{\bold{i},x,y},\bold{j}$\  \  \ \qquad \qquad
$+\sum\limits_{i\gets y\to j}
a_{\overset\curvearrowright\to{ix}}^{(2)}
a_{\overset\curvearrowleft\to{iy}}^{(2)}
=+(t+1)t^{2}$

\smallskip

$\overset+++\to{yxj}:$\ \ \ \ \ $\bold{i},*,
\widehat{x,y,\bold{j}},x,y,\bold{j}$\  \  \ \qquad \qquad
$+\sum\limits_{i\to x\gets j}a_{ix}^{(2)}a_{xj}^{(2)}
+\sum\limits_{i\gets x\gets j}a_{ix}^{(2)}a_{xj}^{(2)}=$

\qquad\quad\qquad\qquad\qquad\qquad\qquad \qquad \quad $+(2t+1)(t+1)t$

\smallskip

$\overset+++\to{yxk}:$\ \ \ \ \ $\bold{i},
k,x,y,\widehat{k,x,y},\bold{j}$\  \  \  \qquad\qquad
$+\sum\limits_{y\to j}
\sum\limits_{i\to k\gets y}
a_{\overset\curvearrowleft\to{ky}}^{(2)}
=+(2t+1)(t+1)t.$

\bigskip

Hence, the  calculations for the case $p=y$ can be summarized as follows:

\smallskip

$III=(4t+1)(2t+1)(t+1)^{2}t-2(4t+1)(t+1)^{2}t-2(2t+1)(t+1)^{2}t$

\smallskip

$-2*2(t+1)^{2}t^{2}+3(2t+1)(t+1)t+3(t+1)t^{2}=(4t+1)(2t+1)(t+1)^{2}t$

\smallskip
$-4(4t+1)(t+1)^{2}t+3(3t+1)(t+1)t=(4t+1)(2t-3)(t+1)^{2}t$

\smallskip
$+3(3t+1)(t+1)t=2t^{2}(t+1)(4t^{2}-t-2).$

\smallskip

All in all, we obtain
$$Corr_{7}(DR_{4t+3},i,j)=I+II+III=$$
$$2t(t+1)(36t^{3}+33t^{2}+3t+1)+2t^{2}(t+1)(2t+1)(4t-3)+2t^{2}(t+1)(4t^{2}-t-2)=$$
$$2t(t+1)\bigl(48t^{3}+30t^{2}-2t+1\bigr).$$

Thus, $c_{8}(DR_{4t+3},ji)=p_{7}(DR_{4t+3},i,j)=
a_{ij}^{(7)}-Corr_{7}(DR_{4t+3},i,j)=$
$$=2t(t+1)\bigl(16t^{4}+28t^{3}+23t^{2}+9t+1\bigr)-2t(t+1)\bigl(48t^{3}+30t^{2}-2t+1\bigr)=$$
$$=2t(t+1)\bigl(16t^{4}-20t^{3}-7t^{2}+11t\bigr)=2(t+1)(t-1)t^{2}(16t^{2}-4t-11).$$

Now we are in a position to formulate the main result of this section.

\smallskip

{\bf Theorem 1.}  {\sl Any arc in a doubly-regular tournament $DR_{4t+3}$
of order $4t+3$ lies on $2(t+1)(t-1)t^{2}(16t^{2}-4t-11)$ cycles of length $8.$}

\medskip

{\bf Corollary 1.} {\sl For any doubly-regular tournament $DR_{n},$ we have}
$$c_{8}(DR_{n})=\frac{(n+1)n(n-1)(n-3)^{2}(n-7)(n^{2}-7n+1)}{256\times 8}.$$

\bigskip

\noindent{\bf \S 3. The number of $\bold{8}$-cycles in $\bold{RLT_{n}}$}

\medskip

In this section, we determine $c_{8}(RLT_{n})$ with the use
of the method first presented in [13]. This method is based on the
self-similar structure of the regular locally transitive tournament:
the one-vertex-deleted subtournament of
$RLT_{n+2}$ is obtained from $RLT_{n}$ by replacing
one of its vertices  with the transitive tournament $TT_{2}$ of order $2.$
This implies that
$$\frac{n-m+2}{n+2}c_{m}(RLT_{n+2})=\frac{n+m}{n}c_{m}(RLT_{n})+
c_{m-1}(RLT_{n},0)+c_{m}^{(2)}(RLT_{n},0),\eqno(1)$$
where $c_{m-1}(RLT_{n},0)$ is the number of $(m-1)$-cycles
starting and ending at $0$ in $RLT_{n}$ and
$c_{m}^{(2)}(RLT_{n},0)$ is the number of closed $m$-walks each of which
is the concatenation of two cycles $\gamma_{1}$ and $\gamma_{2}$
starting and ending at $0$ with $V(\gamma_{1})\cap V(\gamma_{2})=\{0\}$
and $\ell(\gamma_{1})+\ell(\gamma_{2})=m.$

Let $\beta(m)$ be the coefficient of $z^{m-1}$ in the Maclaurin expansion
of $\tan z.$ In other words,  $\beta(m)=\frac{2^{m}(2^{m}-1)|B_{m}|}{m!},$
where $B_{m}$ is the $m$th Bernoulli number.
It is shown in [13] that if we have already proved that both
$c_{m-1}(RLT_{n})$ and $c_{m}^{(2)}(RLT_{n},0)$ are polynomials
of degree $m-1$ and $m-2,$ respectively, in $n,$ namely,
$c_{m-1}(RLT_{n})=\sum\limits_{k=1}^{m-1}\alpha_{k}^{(m-1)}n^{k}$
and
$c_{m}^{(2)}(RLT_{n},0)=\sum\limits_{k=0}^{m-2}\beta_{k}^{(m)}n^{k}$
for each odd $n\ge 1,$ then
$c_{m}(RLT_{n})$ $=\sum\limits_{k=1}^{m}\alpha_{k}^{(m)}n^{k},$
where
$$\alpha_{m}^{(m)}=\frac{1+(-1)^{\frac{m}{2}}\beta(m)}{m2^{m}}$$
and
the other coefficients $\alpha_{1}^{(m)},...,\alpha_{m-1}^{(m)}$
are determined
with the use of the triangular system of $m-1$ linear equations
$$2(k+1-m)\alpha_{k+1}^{(m)}+\sum\limits_{p=k+2}^{m-1}\alpha_{p}^{(m)}2^{p-k-1}\Bigl\{2\binom{p}{k}-m\binom{p-1}{k}\Bigr\}=
f_{k}^{(m)},\eqno(2)$$
where $k$ runs from $0$ and $m-2$ and
$$f_{k}^{(m)}=(m-1)\alpha_{k+1}^{(m-1)}+\beta_{k}^{(m)}+
\frac{m-k-2}{m}2^{-k-1}\binom{m}{k}
\bigl(1+(-1)^{\frac{m}{2}}\beta(m)\bigr).\eqno(3)$$

Consider now the case $m=8$ in detail.
It is well known that
$$\beta(8)=\frac{2^{8}(2^{8}-1)|B_{8}|}{8!}=\frac{2^{8}(2^{8}-1)}{8!\times 30}=
\frac{17}{315}$$
and hence,
$$\alpha_{8}^{(8)}=\frac{332}{8\times 256\times 315}=
\frac{83}{8\times 64\times 315}.$$
Moreover, for $k=0,...,6,$  we have
$$f_{k}^{(8)}=7\alpha_{k+1}^{(7)}+\beta_{k}^{(8)}+(6-k)2^{-k-2}
\binom{8}{k}\frac{83}{315}.$$

It is proved in [13] that
$$c_{7}(RLT_{n})=\frac{(n+1)n(n-1)(n-3)(n-5)(15n^{2}-111n+127)}{128\times 15\times 7}=$$
$$\sum\limits_{k=1}^{7}\alpha_{k}^{(7)}n^{k}=
-\frac{127}{896}n+\frac{383}{1920}n^{2}
+\frac{19}{384}n^{3}-\frac{35}{192}n^{4}+\frac{35}{384}n^{5}
-\frac{11}{640}n^{6}+\frac{1}{896}n^{7}.$$
Moreover, in the next section, we will show that
$$c_{8}^{(2)}(RLT_{n},0)=\frac{(n+1)(n-1)(n-3)(n-5)(85n^{2}-646n+837)}{5760}=$$
$$\sum\limits_{k=0}^{6}\beta_{k}^{(8)}n^{k}=-\frac{279}{128}
+\frac{2731}{960}n+\frac{1055}{1152}n^{2}
-\frac{251}{96}n^{3}+\frac{1439}{1152}n^{4}-\frac{221}{960}n^{5}
+\frac{17}{1152}n^{6}.$$

Hence,
$$f_{0}^{(8)}=7\times\Bigl(-\frac{127}{896}\Bigr)+\Bigl(-\frac{279}{128}\Bigr)+
\frac{6\times 1\times 83}{4 \times 315}=-\frac{18659}{6720};$$
$$f_{1}^{(8)}=7\times\Bigl(\frac{383}{1920}\Bigr)+
\Bigl(\frac{2731}{960}\Bigr)+
\frac{5\times 8\times 83}{8\times 315}=
\frac{224123}{40320};$$
$$f_{2}^{(8)}=7\times\Bigl(\frac{19}{384}\Bigr)+
\Bigl(\frac{1055}{1152}\Bigr)+
\frac{4\times 28\times 83}{16\times 315}=
\frac{8947}{2880};$$
$$f_{3}^{(8)}=7\times\Bigl(-\frac{35}{192}\Bigr)+
\Bigl(-\frac{251}{96}\Bigr)+
\frac{3\times 56 \times 83}{32\times 315}=
-\frac{2407}{960};$$
$$f_{4}^{(8)}=7\times\Bigl(\frac{35}{384}\Bigr)+
\Bigl(\frac{1439}{1152}\Bigr)+
\frac{2\times 70\times 83}{64\times 315}=
\frac{473}{192};$$
$$f_{5}^{(8)}=7\times\Bigl(-\frac{11}{640}\Bigr)+
\Bigl(-\frac{221}{960}\Bigr)+
\frac{1\times 56\times 83}{128\times 315}=
-\frac{271}{1152};$$
$$f_{6}^{(8)}=7\times\Bigl(\frac{1}{896}\Bigr)+
\Bigl(\frac{17}{1152}\Bigr)=\frac{13}{576}.$$

These results of the computation taken together with $(2)$
yield the system of equations
$$\left\{\aligned -2\alpha_{7}^{(8)}=\frac{13}{576}  \\
                 -4\alpha_{6}^{(8)}-12\alpha_{7}^{(8)}=-\frac{271}{1152}    \\
                 -6\alpha_{5}^{(8)}-20\alpha_{6}^{(8)}-200\alpha_{7}^{(8)}=\frac{473}{192}\\
		 -8\alpha_{4}^{(8)}-24\alpha_{5}^{(8)}-160\alpha_{6}^{(8)}-720\alpha_{7}^{(8)}=-\frac{2407}{960}\\
		 -10\alpha_{3}^{(8)}-24\alpha_{4}^{(8)}-112\alpha_{5}^{(8)}-400\alpha_{6}^{(8)}-1248\alpha_{7}^{(8)}=\frac{8947}{2880}\\
		 -12\alpha_{2}^{(8)}-20\alpha_{3}^{(8)}-64\alpha_{4}^{(8)}-176\alpha_{5}^{(8)}-448\alpha_{6}^{(8)}-1088\alpha_{7}^{(8)}=\frac{224123}{40320}\\
		 -14\alpha_{1}^{(8)}-12\alpha_{2}^{(8)}-24\alpha_{3}^{(8)}
-48\alpha_{4}^{(8)}-96\alpha_{5}^{(8)}-192\alpha_{6}^{(8)}
-384\alpha_{7}^{(8)}=-\frac{18659}{6720}.\endaligned\right.$$
Resolving this system yields
$\alpha_{1}^{(8)}=\frac{85}{256},$
$\alpha_{2}^{(8)}=-\frac{6439}{10752},$
$\alpha_{3}^{(8)}=-\frac{13}{576},$
$\alpha_{4}^{(8)}=\frac{11651}{23040},$
$\alpha_{5}^{(8)}=-\frac{791}{2304},$
$\alpha_{6}^{(8)}=\frac{427}{4608},$
$\alpha_{7}^{(8)}=-\frac{13}{1152}.$
Recall also that
$\alpha_{8}^{(8)}=\frac{83}{161280}.$
Factoring the polynomial
$\sum\limits_{k=1}^{8}\alpha_{k}^{(8)}n^{k}$ obtained
yields the following proposition.

{\bf Theorem 2.}  {\sl  For each odd $n\ge 1,$ we have}
$$c_{8}(RLT_{n})=
\frac{(n+1)n(n-1)(n-3)(n-5)(n-7)
(83n^{2}-575n+510)}{161280}.$$

It is not difficult to check that $c_{8}(RLT_{n})-c_{8}(DR_{n})=$
$$\frac{(n+1)n(n-1)(n-3)(n-7)
(17n^{3}-810n^{2}+6610n-9255)}{8\times 80640}.$$
For the cubic polynomial $f(z)=17z^{3}-810z^{2}+6610z-9255,$
we have
$f(1)<0,$ but $f(2)>0,$
$f(8)>0,$ while $f(9)<0,$ and
$f(37)<0,$ but $f(38)>0.$
Hence, the difference $c_{8}(RLT_{n})-c_{8}(DR_{n})$
is positive for each $n\ge 39.$

{\bf Corollary 2.}  {\sl For odd $n\ge 9,$ the strict inequality }
$$c_{8}(DR_{n})<c_{8}(RLT_{n})$$
{\sl holds iff $n\ge 39.$}

Let $A$ be the adjacency matrix of $T$ and $\sigma_{m}(T)$
be the sum of the $m$th powers ($m$th moment)
of the imaginary parts of its (non-Perron) eigenvalues.
(Obviously, $\sigma_{m}(T)=0$ for any odd $m.$)
It was shown in [14] that there exists an absolute
constant $C_{m}$ (not depending on $n$) such that the inequality
$$\Bigl|mc_{m}(T)-\frac{n^{m}}{2^{m}}
-(-1)^{\frac{m}{2}}\sigma_{m}(T) \Bigr|\le C_{m}n^{m-1}$$
holds for any regular tournament $T$ of order $n.$
According to [12], for $m=4,$
the maximum of $\sigma_{m}(T)$
in the class of regular tournaments of order $n$ is attained iff $T=RLT_{n}.$
We conjecture that the same holds for each even $m\ge 4.$
Note that $\sigma_{m}(RLT_{n})=\frac{\beta(m)}{2^{m}}n^{m}+O(n^{m-1}).$
It is not difficult to check that for any even $m\ge 4$ and $T\in {\Cal R}_{n},$
we have $\sigma_{m}(T)\ge \frac{(n-1)n^{\frac{m}{2}}}{2^{m}}$
with equality holding iff $T$ is a doubly-regular tournament of order $n.$
All these arguments and results of [5] allow us to present the following conjecture on $c_{8}(T).$

{\bf Conjecture 1.}  {\sl  For any regular tournament $T$
of sufficiently large order $n,$ we have}
$$c_{8}(DR_{n})\le c_{8}(T)\le c_{8}(RLT_{n}).$$

Corollary 2 shows that
the assumption that $n$ must be large enough is essential.
Our computer processing of B. McKay's files of tournaments
implies that for $n=9$ and $n=11,$ the maximum of $c_{8}(T)$ in a
tournament of order $n$ is achieved for the unique rotational
nearly-doubly-regular tournament $R(2,3,4,8)$
of order $n=9$ and for the quadratic residue tournament $QR_{11},$
which is the unique
doubly-regular tournament of order $n=11,$ respectively.
Moreover, the computer processing of B. McKay's file of
regular tournaments of order $n=13$
implies that the maximum of $c_{8}(T)$ in the class of regular
tournaments of order $n=13$ is attained at some nearly-doubly-regular
tournament, which is not vertex-transitive.
(The interested reader can find its description
in Appendix A.1).
In turn, for $n=11,13,$
the minimum of $c_{8}(T)$
is attained at regular tournaments which are closer,
in their local properties, to $RLT_{n}$ than to $DR_{n}$ or $NDR_{n}$
(see the distribution of $3$-cycles among the out-sets
and in-sets of their vertices in Appendix A.2).
Moreover, for $n=11,13,$ the value of
$c_{8}(RLT_{n})$ is closer to the minimum of $c_{8}(T)$ than to
its maximum.
\footnote[2]{For $n=9,$  the unique minimizer
is isomorphic to the wreath product
$\Delta\circ\Delta,$ which has the largest automorphism group among
all tournaments on $n$ vertices (its order is $81$) and also is  isomorphic to
any Cayley tournament on the abelian group
${\Bbb Z}_{3}\times {\Bbb Z}_{3}$. In turn, the maximum
is attained only at the unique rotational nearly-doubly-regular tournament
$RNDR_{9}$ of order $9.$
Our computer search shows that
$c_{8}(\Delta\circ\Delta)=405\le c_{8}(RLT_{9})=441\le
c_{8}(RNDR_{9})=477.$ Hence,
$c_{8}(RLT_{9})$ is the arithmetical mean of
the corresponding extreme values.}
So, the local properties
of the extreme "points" interchange under the transition from small
to large values of $n.$ Corollary 2 allows us to treat
$n=39$ as the critical value
for the considered extremal combinatorial problem.

\bigskip

\noindent{\bf \S 4. The number $\bold{c_{8}^{(2)}(RLT_{n},0)}$}

\medskip

Recall that the out-set and in-set of each vertex in $RLT_{n}$
induce transitive tournaments.
In the sequel, we assume that
the subtournaments induced by the out-set $N^{+}(0)$ and the in-set
$N^{-}(0)$ in $RLT_{2\delta+1}$
are $TT_{\delta}(1^{+},...,{\delta}^{+})$
and $TT_{\delta}(1^{-},...,{\delta}^{-}),$ respectively,
where $\delta=\frac{n-1}{2}$  and one can assume that
$k^{-}=k^{+}+\delta$ in the ring of residues modulo $2\delta+1.$
Obviously, the condition
$N^{+}(k^{+})\cap N^{-}(0)=\{1^{-},...,k^{-}\}$
uniquely determines the arcs between $N^{+}(0)$ and $N^{-}(0)$.

Let $L_{m}(RLT_{n},0)$ be the set of $m$-cycles
starting and ending at $0$ in $RLT_{n}.$
Obviously, $c_{m}(RLT_{n},0)=|L_{m}(RLT_{n},0)|.$
Denote by $c_{k;h}^{(2)}(RLT_{n},0)$ the number of
concatenations $\gamma_{1}\gamma_{2}$
of $\gamma_{1}\in L_{k}(RLT_{n},0)$ and $\gamma_{2}\in L_{h}(RLT_{n},0)$
with $V(\gamma_{1})\cap V(\gamma_{2})=\{0\}.$
Then
$$c_{8}^{(2)}(RLT_{n},0)=c_{3;5}^{(2)}(RLT_{n},0)+
c_{4;4}^{(2)}(RLT_{n},0)+c_{5;3}^{(2)}(RLT_{n},0).$$
Obviously,
$$c_{k;h}^{(2)}(RLT_{n},0)=c_{k}(RLT_{n},0)\times c_{h}(RLT_{n},0)-
\gamma_{k;h}^{(2)}(RLT_{n},0),$$
where
$\gamma_{k;h}^{(2)}(RLT_{n},0)$ is the number of
concatenations $\gamma_{1}\gamma_{2}$
of $\gamma_{1}\in L_{k}(RLT_{n},0)$ and $\gamma_{2}\in L_{h}(RLT_{n},0)$
with $V(\gamma_{1})\cap V(\gamma_{2})\neq\{0\}.$

To determine $\gamma_{3;5}^{(2)}(RLT_{n},0),$
we consider concatenations
$\gamma_{1}\gamma_{2}$ of $\gamma_{1}\in L_{3}(RLT_{n},$ $0)$
and $\gamma_{2}\in L_{5}(RLT_{n},0)$ with
$V(\gamma_{1})\cap N^{+}(0)\subset V(\gamma_{2}).$
Let us first count the closed $8$-walks of the form
${\bold 0},\star^{+},*^{-},{\bold 0},\star^{+},*^{+},*^{+},*^{-},{\bold 0},$
where $\star^{+}$ (or $\star^{-}$) denotes  coinciding vertices
in $N^{+}(0)\cap V(\gamma_{1})$ and $N^{+}(0)\cap V(\gamma_{2})$
(or in $N^{-}(0)\cap V(\gamma_{1})$ and $N^{-}(0)\cap V(\gamma_{2})$)
and $*^{+}$ (or $*^{-}$) is an arbitrary admissible
vertex in $N^{+}(0)$ (or $N^{-}(0)$). (The word
"admissible" here means that the right (left) neighbour of $*^{+}$
(or $*^{-}$) in the sequence
is its out (in)-neighbour in $RLT_{n}.$)
Replace $\star^{+}$ by $k^{+}.$
Obviously, the number $c_{k}$ of the closed $8$-walks
of the form ${\bold 0},k^{+},*^{-},{\bold 0},k^{+},*^{+},*^{+},*^{-},{\bold 0}$
is equal to the product of
the number $c_{k,1}$ of $3$-circuits of the form ${\bold 0},k^{+},*^{-},{\bold 0}$
and the number $c_{k,2}$ of $5$-circuits of the form
${\bold 0},k^{+},*^{+},*^{+},*^{-},{\bold 0}.$ Recall that
$N^{+}(k^{+})\cap N^{-}(0)=\{1^{-},...,k^{-}\}.$
Hence, $c_{k,1}=k.$ To determine $c_{k,2},$ we consider
$5$-circuits of the form ${\bold 0},k^{+},t^{+},s^{+},*^{-},{\bold 0}.$
Here $t$ varies between $k+1$ to $\delta$ and $s$
runs from $t+1$ to $\delta$ ($*^{-}$ denotes a
out-neighbour of $s^{+}$ in $N^{-}(0)$).
As we have seen above,
the vertex $s^{+}$ has $s$ out-neighbours in $N^{-}(0).$
Therefore, $c_{k,2}=\sum\limits_{t=k+1}^{\delta}
\sum\limits_{s=t+1}^{\delta}s$ and $c_{k}=k\sum\limits_{t=k+1}^{\delta}
\sum\limits_{s=t+1}^{\delta}s.$ This implies that
the number of the closed $8$-walks of the form
${\bold 0},\star^{+},*^{-},{\bold 0},\star^{+},*^{+},*^{+},*^{-},{\bold 0}$
equals $\sum\limits_{k=1}^{\delta}c_{k}=\sum\limits_{k=1}^{\delta}
k\sum\limits_{t=k+1}^{\delta}\sum\limits_{s=t+1}^{\delta}s.$
This case (as well as all the others) can be reflected as follows:

\bigskip

${\bold 0},\star^{+},*^{-},{\bold 0},\star^{+},*^{+},*^{+},*^{-},{\bold 0}=
\bigcup\limits_{k=1}^{\delta}
\bigcup\limits_{t=k+1}^{\delta}\bigcup\limits_{s=t+1}^{\delta}
{\bold 0},{k}^{+},*^{-},{\bold 0},{k}^{+},t^{+},s^{+},*^{-},{\bold 0}$   $(a1)$

\smallskip

$\sum\limits_{k=1}^{\delta}k
\sum\limits_{t=k+1}^{\delta}
\sum\limits_{s=t+1}^{\delta}s=
\sum\limits_{k=1}^{\delta}
\frac{k(\delta-1-k)(\delta-k)(2\delta+2+k)}{6}=
\frac{(\delta-2)(\delta-1)\delta(\delta+1)(4\delta+3)}{120}$

\medskip

${\bold 0},\star^{+},*^{-},{\bold 0},*^{+},\star^{+},*^{+},*^{-},{\bold 0}
=\bigcup\limits_{k=1}^{\delta}
\bigcup\limits_{t=k+1}^{\delta}
{\bold 0},k^{+},*^{-},{\bold 0},*^{+},k^{+},t^{+},*^{-},{\bold 0}$  \qquad\ $(a2)$

\smallskip

$\sum\limits_{k=1}^{\delta}k(k-1)\sum\limits_{t=k+1}^{\delta}t
=\sum\limits_{k=1}^{\delta}\frac{k(k-1)(\delta+1+k)(\delta-k)}{2}=
\frac{(\delta-2)(\delta-1)\delta(\delta+1)(4\delta+3)}{60}$

\medskip

${\bold 0},\star^{+},*^{-},{\bold 0},*^{+},*^{+},\star^{+},*^{-},{\bold 0}=
\bigcup\limits_{k=1}^{\delta}\bigcup\limits_{s=1}^{k-1}
{\bold 0},k^{+},*^{-},{\bold 0},*^{+},s^{+},k^{+},*^{-},{\bold 0}$ \qquad \quad \  $(a3)$

\smallskip

$\sum\limits_{k=1}^{\delta}k\sum\limits_{s=1}^{k-1}(s-1)k=
\sum\limits_{k=1}^{\delta}\frac{k^{2}(k-1)(k-2)}{2}=
\frac{(\delta-2)(\delta-1)\delta(\delta+1)(4\delta+3)}{40}$

\medskip

${\bold 0},\star^{+},*^{-},{\bold 0},\star^{+},*^{+},*^{-},*^{-},{\bold 0}
=\bigcup\limits_{k=1}^{\delta}
\bigcup\limits_{t=k+1}^{\delta}\bigcup\limits_{s=1}^{t}
{\bold 0},k^{+},*^{-},{\bold 0},k^{+},t^{+},s^{-},*^{-},{\bold 0}$  \quad  $(a4)$

\smallskip

$\sum\limits_{k=1}^{\delta}k
\sum\limits_{t=k+1}^{\delta}
\sum\limits_{s=1}^{t}(\delta-s)=
\sum\limits_{k=1}^{\delta}\frac{k(\delta-k)(2\delta^{2}+2\delta k-2-3k-k^{2})}{6}=
\frac{(\delta-1)\delta(\delta+1)(9\delta^{2}-5\delta-6)}{120}$

\medskip

${\bold 0},\star^{+},*^{-},{\bold 0},*^{+},\star^{+},*^{-},*^{-},{\bold 0}
=\bigcup\limits_{k=1}^{\delta}
\bigcup\limits_{s=1}^{k}{\bold 0},k^{+},*^{-},{\bold 0},*^{+},k^{+},s^{-},*^{-},{\bold 0}$ \qquad \quad \ $(a5)$

\smallskip

$\sum\limits_{k=1}^{\delta}k(k-1)
\sum\limits_{s=1}^{k}(\delta-s)=
\sum\limits_{k=1}^{\delta}\frac{k^{2}(k-1)(2\delta-1-k)}{2}=\frac{(\delta-1)\delta(\delta+1)(9\delta^{2}-5\delta-6)}{60}$

\medskip

${\bold 0},\star^{+},*^{-},{\bold 0},\star^{+},*^{-},*^{+},*^{-},{\bold 0}=
\bigcup\limits_{k=1}^{\delta}
\bigcup\limits_{s=1}^{k}\bigcup\limits_{t=1}^{s-1}
{\bold 0},k^{+},*^{-},{\bold 0},k^{+},s^{-},t^{+},*^{-},{\bold 0}$ \quad \ \ $(a6)$

\smallskip

$\sum\limits_{k=1}^{\delta}k\sum\limits_{s=1}^{k}\sum\limits_{t=1}^{s-1}t=
\sum\limits_{k=1}^{\delta}\frac{k^{2}(k+1)(k-1)}{6}=
\frac{(\delta-1)\delta(\delta+1)(\delta+2)(2\delta+1)}{60}$

\medskip

${\bold 0},\star^{+},*^{-},{\bold 0},*^{+},*^{-},\star^{+},*^{-},{\bold 0}
=\bigcup\limits_{k=1}^{\delta}\bigcup\limits_{s=k+1}^{\delta}
{\bold 0},k^{+},*^{-},{\bold 0},*^{+},s^{-},k^{+},*^{-},{\bold 0}$ \qquad \ $(a7)$

\smallskip

$\sum\limits_{k=1}^{\delta}k
\sum\limits_{s=k+1}^{\delta}(\delta-s+1)k=\sum\limits_{k=1}^{\delta}
\frac{k^{2}(\delta-k)(\delta-k+1)}{2}=
\frac{(\delta-1)\delta(\delta+1)(\delta+2)(2\delta+1)}{120}$

\medskip

${\bold 0},\star^{+},*^{-},{\bold 0},\star^{+},*^{-},*^{-},*^{-},{\bold 0}
=\bigcup\limits_{k=1}^{\delta}
\bigcup\limits_{s=1}^{k}\bigcup\limits_{t=s+1}^{\delta}
{\bold 0},k^{+},*^{-},{\bold 0},k^{+},s^{-},t^{-},*^{-},{\bold 0}$  \ \ \  $(a8)$

\smallskip

$\sum\limits_{k=1}^{\delta}k
\sum\limits_{s=1}^{k}
\sum\limits_{t=s+1}^{\delta}(\delta-t)=
\sum\limits_{k=1}^{\delta}
\frac{k^{2}(3\delta^{2}-6\delta+2-3\delta k +3k+k^{2})}{6}=
\frac{(\delta-1)\delta(\delta+1)(\delta-2)(3\delta+1)}{40}.$

\medskip

Since any tournament admits no $2$-cycles, for any two $i,j=1,...,8,$
the subsets $(ai)$ and $(aj)$ contain no common
closed $8$-walks. Summing from $(a1)$ to $(a8)$ yields the sum
$$A=\frac{(\delta-1)\delta(\delta+1)(\delta+2)(2\delta+1)}{40}+
\frac{(\delta-1)\delta(\delta+1)(\delta-2)(3\delta+1)}{40}+$$
$$\frac{(\delta-2)(\delta-1)\delta(\delta+1)(4\delta+3)}{20}+
\frac{(\delta-1)\delta(\delta+1)(9\delta^{2}-5\delta-6)}{40}=$$
$$\frac{(\delta-1)\delta(\delta+1)(22\delta^{2}-15\delta-18)}{40}.$$

By duality,
the number of concatenations
$\gamma_{1}\gamma_{2}$ of $\gamma_{1}\in L_{3}(RLT_{n},0)$
and $\gamma_{2}\in L_{5}(RLT_{n},$ $0)$ with
$V(\gamma_{1})\cap N^{-}(0)\subset V(\gamma_{2})$
equals
the number of concatenations
$\gamma_{1}\gamma_{2}$ of $\gamma_{1}\in L_{3}(RLT_{n},0)$
and $\gamma_{2}\in L_{5}(RLT_{n},0)$ with
$V(\gamma_{1})\cap N^{+}(0)\subset V(\gamma_{2}).$
However, concatenations $\gamma_{1}\gamma_{2}$ of
$\gamma_{1}\in L_{3}(RLT_{n},0)$ and $\gamma_{2}\in L_{5}(RLT_{n},0)$
with $V(\gamma_{1})\subset V(\gamma_{2})$
are counted twice.
Let us determine the number
of such walks.

\medskip

${\bold 0},\star^{+},\star^{-},{\bold 0},\star^{+},*^{+},*^{+},\star^{-},{\bold 0}=
\bigcup\limits_{k=1}^{\delta}
\bigcup\limits_{s=1}^{k}\bigcup\limits_{t=k+1}^{\delta}
{\bold 0},k^{+},s^{-},{\bold 0},k^{+},t^{+},*^{+},s^{-},{\bold 0}$ \quad \ $(b1)$

\smallskip

$\sum\limits_{k=1}^{\delta}\sum\limits_{s=1}^{k}\sum\limits_{t=k+1}^{\delta}
(\delta-t)=\sum\limits_{k=1}^{\delta}\frac{k(\delta-k-1)(\delta-k)}{2}=
\frac{(\delta-2)(\delta-1)\delta(\delta+1)}{24}$

\medskip

${\bold 0},\star^{+},\star^{-},{\bold 0},*^{+},\star^{+},*^{+},\star^{-},{\bold 0}=
\bigcup\limits_{k=1}^{\delta}
\bigcup\limits_{s=1}^{k}
{\bold 0},k^{+},s^{-},{\bold 0},*^{+},k^{+},*^{+},s^{-},{\bold 0}$ \quad\quad \quad \ $(b2)$

\smallskip

$\sum\limits_{k=1}^{\delta}\sum\limits_{s=1}^{k}(k-1)(\delta-k)
=\sum\limits_{k=1}^{\delta}k(k-1)(\delta-k)=
\frac{(\delta-2)(\delta-1)\delta(\delta+1)}{12}$

\medskip

${\bold 0},\star^{+},\star^{-},{\bold 0},*^{+},*^{+},\star^{+},\star^{-},{\bold 0}=
\bigcup\limits_{k=1}^{\delta}\bigcup\limits_{s=1}^{k}\bigcup\limits_{t=1}^{k-1}
{\bold 0},k^{+},s^{-},{\bold 0},t^{+},*^{+},k^{+},s^{-},{\bold 0}$ \qquad \ $(b3)$

\smallskip

$\sum\limits_{k=1}^{\delta}\sum\limits_{s=1}^{k}\sum\limits_{t=1}^{k-1}
(k-1-t)
=\sum\limits_{k=1}^{\delta}\frac{k(k-1)(k-2)}{2}=
\frac{(\delta-2)(\delta-1)\delta(\delta+1)}{8}$

\medskip

${\bold 0},\star^{+},\star^{-},{\bold 0},\star^{+},*^{+},\star^{-},*^{-},{\bold 0}=
\bigcup\limits_{k=1}^{\delta}\bigcup\limits_{s=1}^{k}
{\bold 0},k^{+},s^{-},{\bold 0},k^{+},*^{+},s^{-},*^{-},{\bold 0}$ \qquad \quad \ $(b4)$

\smallskip

$\sum\limits_{k=1}^{\delta}\sum\limits_{s=1}^{k}
(\delta-k)(\delta-s)
=\sum\limits_{k=1}^{\delta}\frac{k(\delta-k)(2\delta-1-k)}{2}=
\frac{(\delta-1)\delta(\delta+1)(3\delta-2)}{24}$

\medskip

${\bold 0},\star^{+},\star^{-},{\bold 0},\star^{+},*^{+},*^{-},\star^{-},{\bold 0}=
\bigcup\limits_{k=1}^{\delta}\bigcup\limits_{s=1}^{k}
{\bold 0},k^{+},s^{-},{\bold 0},k^{+},*^{+},*^{-},s^{-},{\bold 0}$ \qquad\quad \ $(b5)$

\smallskip

$\sum\limits_{k=1}^{\delta}\sum\limits_{s=1}^{k}
(\delta-k)(s-1)
=\sum\limits_{k=1}^{\delta}\frac{(\delta-k)k(k-1)}{2}=
\frac{(\delta-2)(\delta-1)\delta(\delta+1)}{24}$

\medskip

${\bold 0},\star^{+},\star^{-},{\bold 0},*^{+},\star^{+},\star^{-},*^{-},{\bold 0}=
\bigcup\limits_{k=1}^{\delta}
\bigcup\limits_{s=1}^{k}
{\bold 0},k^{+},s^{-},{\bold 0},*^{+},k^{+},s^{-},*^{-},{\bold 0}$ \qquad\quad \ $(b6)$

\smallskip

$\sum\limits_{k=1}^{\delta}\sum\limits_{s=1}^{k}
(k-1)(\delta-s)
=\sum\limits_{k=1}^{\delta}
\frac{(k-1)k(2\delta-1-k)}{2}=
\frac{(\delta-1)\delta(\delta+1)(5\delta-6)}{24}$

\medskip

${\bold 0},\star^{+},\star^{-},{\bold 0},*^{+},\star^{+},*^{-},\star^{-},{\bold 0}
=\bigcup\limits_{k=1}^{\delta}
\bigcup\limits_{s=1}^{k}
{\bold 0},k^{+},s^{-},{\bold 0},*^{+},k^{+},*^{-},s^{-},{\bold 0}$ \qquad\quad \ $(b7)$

\smallskip

Reversing $\gamma_{2}$ implies that $|(b7)|=|(b4)|=\frac{(\delta-1)\delta(\delta+1)(3\delta-2)}{24}$

\medskip

${\bold 0},\star^{+},\star^{-},{\bold 0},\star^{+},\star^{-},*^{+},*^{-},{\bold 0}=
\bigcup\limits_{k=1}^{\delta}
\bigcup\limits_{s=1}^{k}\bigcup\limits_{t=1}^{s-1}
{\bold 0},k^{+},s^{-},{\bold 0},k^{+},s^{-},t^{+},*^{-},{\bold 0}$ \qquad \ $(b8)$

\smallskip

$\sum\limits_{k=1}^{\delta}\sum\limits_{s=1}^{k}
\sum\limits_{t=1}^{s-1}t
=\sum\limits_{k=1}^{\delta}
\sum\limits_{s=1}^{k}\frac{s(s-1)}{2}=
\sum\limits_{k=1}^{\delta}\frac{(k-1)k(k+1)}{6}
=\frac{(\delta-1)\delta(\delta+1)(\delta+2)}{24}$

\medskip

${\bold 0},\star^{+},\star^{-},{\bold 0},\star^{+},*^{-},*^{+},\star^{-},{\bold 0}
=\bigcup\limits_{k=1}^{\delta}
\bigcup\limits_{s=1}^{k}\bigcup\limits_{t=s+1}^{k}
{\bold 0},k^{+},s^{-},{\bold 0},k^{+},t^{-},*^{+},s^{-},{\bold 0}$ \quad \   $(b9)$

\smallskip

$\sum\limits_{k=1}^{\delta}\sum\limits_{s=1}^{k}
\sum\limits_{t=s+1}^{k}(t-s)=
\sum\limits_{k=1}^{\delta}\frac{(k-1)k(k+1)}{6}=\frac{(\delta-1)\delta(\delta+1)(\delta+2)}{24}$

\medskip

${\bold 0},\star^{+},\star^{-},{\bold 0},*^{+},*^{-},\star^{+},\star^{-},{\bold 0}=
\bigcup\limits_{k=1}^{\delta}
\bigcup\limits_{s=1}^{k}\bigcup\limits_{t=k+1}^{\delta}
{\bold 0},k^{+},s^{-},{\bold 0},*^{+},t^{-},k^{+},s^{-},{\bold 0}$  \ \ $(b10)$

\smallskip

Reversing $\gamma_{2}$ implies that $|(b10)|=|(b8)|=\frac{(\delta-1)\delta(\delta+1)(\delta+2)}{24}$

\medskip

${\bold 0},\star^{+},\star^{-},{\bold 0},\star^{+},\star^{-},*^{-},*^{-},{\bold 0}=
\bigcup\limits_{k=1}^{\delta}
\bigcup\limits_{s=1}^{k}\bigcup\limits_{t=s+1}^{\delta}
{\bold 0},k^{+},s^{-},{\bold 0},k^{+},s^{-},t^{-},*^{-},{\bold 0}$  \ \  $(b11)$

\smallskip

Reversing $\gamma_{2}$ implies that $|(b11)|=|(b3)|=\frac{(\delta-2)(\delta-1)\delta(\delta+1)}{8}$

\medskip

${\bold 0},\star^{+},\star^{-},{\bold 0},\star^{+},*^{-},\star^{-},*^{-},{\bold 0}=
\bigcup\limits_{k=1}^{\delta}
\bigcup\limits_{s=1}^{k}
{\bold 0},k^{+},s^{-},{\bold 0},k^{+},*^{-},s^{-},*^{-},{\bold 0}$ \qquad\quad    $(b12)$

\smallskip

Reversing $\gamma_{2}$ implies that $|(b12)|=|(b2)|=\frac{(\delta-2)(\delta-1)\delta(\delta+1)}{12}$

\medskip

${\bold 0},\star^{+},\star^{-},{\bold 0},\star^{+},*^{-},*^{-},\star^{-},{\bold 0}=
\bigcup\limits_{k=1}^{\delta}
\bigcup\limits_{s=1}^{k}\bigcup\limits_{t=1}^{s-1}
{\bold 0},k^{+},s^{-},{\bold 0},k^{+},t^{-},*^{-},s^{-},{\bold 0}$ \quad \ \  \ $(b13)$

\smallskip

Reversing $\gamma_{2}$ implies that $|(b13)|=|(b1)|=\frac{(\delta-2)(\delta-1)\delta(\delta+1)}{24}.$

\medskip

Summing from $(b1)$ to $(b13)$ yields the sum
$$B=\frac{13(\delta-2)(\delta-1)\delta(\delta+1)}{24}+
\frac{(\delta-1)\delta(\delta+1)(\delta+2)}{8}+$$
$$\frac{(\delta-1)\delta(\delta+1)(3\delta-2)}{12}+
\frac{(\delta-1)\delta(\delta+1)(5\delta-6)}{24}=
\frac{(\delta-1)\delta(\delta+1)(9\delta-10)}{8}.$$

Thus, $c_{3;5}^{(2)}(RLT_{n},0)=c_{3}(RLT_{n},0)\times c_{5}(RLT_{n},0)
-2A+B=$
$$\frac{\delta(\delta+1)}{2}\times\frac{\delta(\delta-1)(\delta+1)
(3\delta-4)}{6}-
\frac{(\delta-1)\delta(\delta+1)(22\delta^{2}-15\delta-18)}{20}$$
$$+\frac{(\delta-1)\delta(\delta+1)(9\delta-10)}{8}=
\frac{(\delta-2)(\delta-1)\delta(\delta+1)(30\delta^{2}-82\delta+21)}{120}.$$

Let us determine now $\gamma_{4;4}^{(2)}(RLT_{n},0).$
We start with concatenations
$\gamma_{1}\gamma_{2}$ of $\gamma_{1},\gamma_{2}$ $\in L_{4}(RLT_{n},0)$
such that $V(\gamma_{1})\cap N^{+}(0)\cap V(\gamma_{2})\neq \emptyset.$

\bigskip

${\bold 0},\star^{+},*^{+},*^{-},{\bold 0},\star^{+},*^{+},*^{-},{\bold 0}=
\bigcup\limits_{k=1}^{\delta}
\bigcup\limits_{t=k+1}^{\delta}\bigcup\limits_{s=k+1}^{\delta}
{\bold 0},k^{+},t^{+},*^{-},{\bold 0},k^{+},s^{+},*^{-},{\bold 0}$  $(a1)$

\smallskip

$\sum\limits_{k=1}^{\delta}\sum\limits_{t=k+1}^{\delta}
\sum\limits_{s=k+1}^{\delta}ts=
\sum\limits_{k=1}^{\delta}\frac{(\delta+1+k)^{2}(\delta-k)^{2}}{4}=
\frac{(\delta-1)\delta(\delta+1)(8\delta^{2}+5\delta-2)}{60}$

\medskip

${\bold 0},\star^{+},*^{+},*^{-},{\bold 0},*^{+},\star^{+},*^{-},{\bold 0}=
\bigcup\limits_{k=1}^{\delta}
\bigcup\limits_{t=k+1}^{\delta}
{\bold 0},k^{+},t^{+},*^{-},{\bold 0},*^{+},k^{+},*^{-},{\bold 0}$ \qquad\ \ \ $(a2)$

\smallskip

$\sum\limits_{k=1}^{\delta}\sum\limits_{t=k+1}^{\delta}t(k-1)k=
\sum\limits_{k=1}^{\delta}\frac{(k-1)k(\delta+1+k)(\delta-k)}{2}=
\frac{(\delta-2)(\delta-1)\delta(\delta+1)(4\delta+3)}{60}$

\medskip

${\bold 0},*^{+},\star^{+},*^{-},{\bold 0},\star^{+},*^{+},*^{-},{\bold 0}=
\bigcup\limits_{k=1}^{\delta}
\bigcup\limits_{t=k+1}^{\delta}
{\bold 0},*^{+},k^{+},*^{-},{\bold 0},k^{+},t^{+},*^{-},{\bold 0}$ \qquad\ \ \ $(a3)$

\smallskip

Interchanging $\gamma_{1}$ and $\gamma_{2}$ implies that
$|(a3)|=|(a2)|=\frac{(\delta-2)(\delta-1)\delta(\delta+1)(4\delta+3)}{60}$

\medskip

${\bold 0},*^{+},\star^{+},*^{-},{\bold 0},*^{+},\star^{+},*^{-},{\bold 0}=
\bigcup\limits_{k=1}^{\delta}
{\bold 0},*^{+},k^{+},*^{-},{\bold 0},*^{+},k^{+},*^{-},{\bold 0}$ \qquad\qquad\quad  $(a4)$

\smallskip

$\sum\limits_{k=1}^{\delta}(k-1)k(k-1)k=
\sum\limits_{k=1}^{\delta}k^{2}(k-1)^{2}=
\frac{(\delta-1)\delta(\delta+1)(3\delta^{2}-2)}{15}$

\medskip

${\bold 0},\star^{+},*^{+},*^{-},{\bold 0},\star^{+},*^{-},*^{-},{\bold 0}=
\bigcup\limits_{k=1}^{\delta}
\bigcup\limits_{t=k+1}^{\delta}\bigcup\limits_{s=1}^{k}
{\bold 0},k^{+},t^{+},*^{-},{\bold 0},k^{+},s^{-},*^{-},{\bold 0}$  \quad   $(a5)$

\smallskip

$\sum\limits_{k=1}^{\delta}\sum\limits_{t=k+1}^{\delta}
\sum\limits_{s=1}^{k}t(\delta-s)=\sum\limits_{k=1}^{\delta}
\frac{(\delta+1+k)(\delta-k)k(2\delta-1-k)}{4}=\frac{(\delta-1)\delta(\delta+1)(11\delta^{2}-4)}{120}$

\medskip

${\bold 0},*^{+},\star^{+},*^{-},{\bold 0},\star^{+},*^{-},*^{-},{\bold 0}=
\bigcup\limits_{k=1}^{\delta}\bigcup\limits_{s=1}^{k}
{\bold 0},*^{+},k^{+},*^{-},{\bold 0},k^{+},s^{-},*^{-},{\bold 0}$ \qquad\quad \ \  $(a6)$

\smallskip

$\sum\limits_{k=1}^{\delta}
\sum\limits_{s=1}^{k}(k-1)k(\delta-s)=
\sum\limits_{k=1}^{\delta}\frac{(k-1)k^{2}(2\delta-1-k)}{2}=\frac{(\delta-1)\delta(\delta+1)(9\delta^{2}-5\delta-6)}{60}$

\medskip

${\bold 0},\star^{+},*^{-},*^{-},{\bold 0},\star^{+},*^{+},*^{-},{\bold 0}=
\bigcup\limits_{k=1}^{\delta}
\bigcup\limits_{s=1}^{k}\bigcup\limits_{t=k+1}^{\delta}
{\bold 0},k^{+},s^{-},*^{-},{\bold 0},k^{+},t^{+},*^{-},{\bold 0}$   \quad $(a7)$

\smallskip

Interchanging $\gamma_{1}$ and $\gamma_{2}$ implies that $|(a7)|=|(a5)|=\frac{(\delta-1)\delta(\delta+1)(11\delta^{2}-4)}{120}$

\medskip

${\bold 0},\star^{+},*^{-},*^{-},{\bold 0},*^{+},\star^{+},*^{-},{\bold 0}=
\bigcup\limits_{k=1}^{\delta}
\bigcup\limits_{s=1}^{k}
{\bold 0},k^{+},s^{-},*^{-},{\bold 0},*^{+},k^{+},*^{-},{\bold 0}$ \qquad\quad \ \ $(a8)$

\smallskip

Interchanging $\gamma_{1}$ and $\gamma_{2}$ implies that $|(a8)|=|(a6)|=\frac{(\delta-1)\delta(\delta+1)(9\delta^{2}-5\delta-6)}{60}$

\medskip

${\bold 0},\star^{+},*^{-},*^{-},{\bold 0},\star^{+},*^{-},*^{-},{\bold 0}=
\bigcup\limits_{k=1}^{\delta}
\bigcup\limits_{s=1}^{k}\bigcup\limits_{t=1}^{k}
{\bold 0},k^{+},s^{-},*^{-},{\bold 0},k^{+},t^{-},*^{-},{\bold 0}$  \qquad \ $(a9)$

\smallskip

$\sum\limits_{k=1}^{\delta}
\sum\limits_{s=1}^{k}\sum\limits_{t=1}^{k}(\delta-s)(\delta-t)=
\sum\limits_{k=1}^{\delta}\frac{k^{2}(2\delta-1-k)^{2}}{4}=\frac{(\delta-1)\delta(\delta+1)(8\delta^{2}-5\delta-2)}{60}.$

\medskip

Note that the subsets $(a1)$ and $(a4)$ admit common
closed $8$-walks. They have

the form ${\bold 0},\star^{+},\star^{+},*^{-},{\bold 0},\star^{+},\star^{+},*^{-},{\bold 0}.$
Obviously,

${\bold 0},\star^{+},\star^{+},*^{-},{\bold 0},\star^{+},\star^{+},*^{-},{\bold 0}=
\bigcup\limits_{k=1}^{\delta}\bigcup\limits_{t=k+1}^{\delta}
{\bold 0},k^{+},t^{+},*^{-},{\bold 0},k^{+},t^{+},*^{-},{\bold 0}.$

Hence, the number of them is equal to $a=
\sum\limits_{k=1}^{\delta}
\sum\limits_{t=k+1}^{\delta}t^{2}=$
$$=\sum\limits_{k=1}^{\delta}\bigl(\frac{\delta(\delta+1)(2\delta+1)}{6}
-\frac{k(k+1)(2k+1)}{6}\bigr)=
\frac{(\delta-1)\delta(\delta+1)(3\delta+2)}{12}.$$

Summing from $(a1)$ to $(a9)$ and then subtracting $a$
yields the following sum
$$A=\frac{(\delta-2)(\delta-1)\delta(\delta+1)(4\delta+3)}{30}+
\frac{(\delta-1)\delta(\delta+1)(3\delta^{2}-2)}{15}+$$
$$\frac{(\delta-1)\delta(\delta+1)(8\delta^{2}-5\delta-2)}{60}
+\frac{(\delta-1)\delta(\delta+1)(8\delta^{2}+5\delta-2)}{60}+$$
$$\frac{(\delta-1)\delta(\delta+1)(9\delta^{2}-5\delta-6)}{30}
+\frac{(\delta-1)\delta(\delta+1)(11\delta^{2}-4)}{60}+$$
$$-\frac{(\delta-1)\delta(\delta+1)(3\delta+2)}{12}=
\frac{(\delta-1)\delta(\delta+1)(13\delta^{2}-7\delta-10)}{12}.$$

By duality,
the number of concatenations
$\gamma_{1}\gamma_{2}$ of $\gamma_{1},\gamma_{2}\in L_{4}(RLT_{n},0)$
such that
$V(\gamma_{1})\cap N^{-}(0)\cap V(\gamma_{2})\neq \emptyset$
equals
the number of concatenations
$\gamma_{1}\gamma_{2}$ of $\gamma_{1},\gamma_{2}\in L_{4}(RLT_{n},0)$
such that
$V(\gamma_{1})\cap N^{+}(0)\cap V(\gamma_{2})\neq \emptyset.$
Note that concatenations $\gamma_{1}\gamma_{2}$ of
$\gamma_{1},\gamma_{2}\in L_{4}(RLT_{n},0)$ having
common vertices both in $N^{+}(0)$ and $N^{-}(0)$ are counted twice.
Let us first determine the number
of such walks with $|V(\gamma_{1})\cap N^{+}(0)|=2.$

\medskip

${\bold 0},\star^{+},*^{+},\star^{-},{\bold 0},\star^{+},*^{+},\star^{-},{\bold 0}=
\bigcup\limits_{k=1}^{\delta}\bigcup\limits_{s=1}^{\delta}
{\bold 0},k^{+},*^{+},s^{-},{\bold 0},k^{+},*^{+},s^{-},{\bold 0}$ \quad \qquad \ \ \ $(b1)$

\smallskip

$\sum\limits_{k=1}^{\delta}
\Bigl(\sum\limits_{s=1}^{k}(\delta-k)^{2}+\sum\limits_{s=k+1}^{\delta}
(\delta-s+1)^{2}\Bigr)=\sum\limits_{k=1}^{\delta}
\frac{(\delta-k)(1+3\delta+2\delta^{2}-3k+2\delta k-4k^{2})}{6}=
\frac{\delta^{2}(\delta^{2}-1)}{6}$

\medskip

${\bold 0},\star^{+},*^{+},\star^{-},{\bold 0},*^{+},\star^{+},\star^{-},{\bold 0}=
\bigcup\limits_{k=1}^{\delta}\bigcup\limits_{s=1}^{k}
{\bold 0},k^{+},*^{+},s^{-},{\bold 0},*^{+},k^{+},s^{-},{\bold 0}$ \qquad\quad \ $(b2)$

\smallskip

$\sum\limits_{k=1}^{\delta}
\sum\limits_{s=1}^{k}(\delta-k)(k-1)=
\sum\limits_{k=1}^{\delta}k(\delta-k)(k-1)=
\frac{(\delta-2)(\delta-1)\delta(\delta+1)}{12}$

\medskip

${\bold 0},*^{+},\star^{+},\star^{-},{\bold 0},\star^{+},*^{+},\star^{-},{\bold 0}=
\bigcup\limits_{k=1}^{\delta}\bigcup\limits_{s=1}^{k}
{\bold 0},*^{+},k^{+},s^{-},{\bold 0},k^{+},*^{+},s^{-},{\bold 0}$ \qquad\quad \ $(b3)$

\smallskip

Interchanging $\gamma_{1}$ and $\gamma_{2}$ implies that $|(b3)|=|(b2)|=\frac{(\delta-2)(\delta-1)\delta(\delta+1)}{12}$

\medskip

${\bold 0},*^{+},\star^{+},\star^{-},{\bold 0},*^{+},\star^{+},\star^{-},{\bold 0}=
\bigcup\limits_{k=1}^{\delta}\bigcup\limits_{s=1}^{k}
{\bold 0},*^{+},k^{+},s^{-},{\bold 0},*^{+},k^{+},s^{-},{\bold 0}$ \qquad\quad \ $(b4)$

\smallskip

$\sum\limits_{k=1}^{\delta}
\sum\limits_{s=1}^{k}(k-1)(k-1)=
\sum\limits_{k=1}^{\delta}k(k-1)^{2}=
\frac{(\delta-1)\delta(\delta+1)(3\delta-2)}{12}$

\medskip

${\bold 0},\star^{+},*^{+},\star^{-},{\bold 0},\star^{+},\star^{-},*^{-},{\bold 0}=
\bigcup\limits_{k=1}^{\delta}\bigcup\limits_{s=1}^{k}
{\bold 0},k^{+},*^{+},s^{-},{\bold 0},k^{+},s^{-},*^{-},{\bold 0}$ \qquad\quad \ $(b5)$

\smallskip

$\sum\limits_{k=1}^{\delta}
\sum\limits_{s=1}^{k}(\delta-k)(\delta-s)=
\sum\limits_{k=1}^{\delta}\frac{(\delta-k)k(2\delta-1-k)}{2}=
\frac{(\delta-1)\delta(\delta+1)(3\delta-2)}{24}$

\medskip

${\bold 0},\star^{+},*^{+},\star^{-},{\bold 0},\star^{+},*^{-},\star^{-},{\bold 0}=
\bigcup\limits_{k=1}^{\delta}\bigcup\limits_{s=1}^{\delta}
{\bold 0},k^{+},*^{+},s^{-},{\bold 0},k^{+},*^{-},s^{-},{\bold 0}$ \quad \qquad \ $(b6)$

\smallskip

$\sum\limits_{k=1}^{\delta}
\Bigl(\sum\limits_{s=1}^{k}(\delta-k)(s-1)
+\sum\limits_{s=k+1}^{\delta}(\delta-s+1)k\Bigr)=
\sum\limits_{k=1}^{\delta}\frac{\delta(\delta-k)k}{2}=
\frac{(\delta-1)\delta^{2}(\delta+1)}{12}$

\medskip

${\bold 0},*^{+},\star^{+},\star^{-},{\bold 0},\star^{+},\star^{-},*^{-},{\bold 0}=
\bigcup\limits_{k=1}^{\delta}\bigcup\limits_{s=1}^{k}
{\bold 0},*^{+},k^{+},s^{-},{\bold 0},k^{+},s^{-},*^{-},{\bold 0}$ \qquad\quad \  $(b7)$

\smallskip

$\sum\limits_{k=1}^{\delta}
\sum\limits_{s=1}^{k}(k-1)(\delta-s)=
\sum\limits_{k=1}^{\delta}\frac{(k-1)k(2\delta-1-k)}{2}=
\frac{(\delta-1)\delta(\delta+1)(5\delta-6)}{24}$

\medskip

${\bold 0},*^{+},\star^{+},\star^{-},{\bold 0},\star^{+},*^{-},\star^{-},{\bold 0}=
\bigcup\limits_{k=1}^{\delta}\bigcup\limits_{s=1}^{k}
{\bold 0},*^{+},k^{+},s^{-},{\bold 0},k^{+},*^{-},s^{-},{\bold 0}$ \qquad\quad \ $(b8)$

\smallskip

Reversing the whole $\gamma_{1}\gamma_{2}$ implies that $|(b8)|=|(b5)|=\frac{(\delta-1)\delta(\delta+1)(3\delta-2)}{24}.$

\medskip

Summing from $(b1)$ to $(b8)$ yields the following sum
$$B=\frac{(\delta-2)(\delta-1)\delta(\delta+1)}{6}+
\frac{(\delta-1)\delta^{2}(\delta+1)}{4}$$
$$\frac{(\delta-1)\delta(\delta+1)(3\delta-2)}{6}
+\frac{(\delta-1)\delta(\delta+1)(5\delta-6)}{24}=
\frac{(\delta-1)\delta(\delta+1)(27\delta-22)}{24}.$$

By duality, the number of concatenations $\gamma_{1}\gamma_{2}$ of
$\gamma_{1},\gamma_{2}\in L_{4}(RLT_{n},0)$
such that $V(\gamma_{1}) \cap N^{+}(0)\cap V(\gamma_{2})\neq \emptyset,$
$V(\gamma_{1})\cap N^{-}(0)\cap V(\gamma_{2})\neq \emptyset,$
and $|V(\gamma_{1})\cap N^{-}(0)|=2$ also equals $B.$
However, the closed $8$-walks of the form
${\bold 0},\star^{+},\star^{+},\star^{-},{\bold 0},\star^{+},\star^{+},
\star^{-},{\bold 0}$ or
${\bold 0},\star^{+},\star^{-},$ $\star^{-},{\bold 0},$ $\star^{+},\star^{-},
\star^{-},{\bold 0}$ are counted twice and hence, should be subtracted
from the sum. Obviously, the number of them
coincides with $c_{4}(RLT_{n},0)$  and hence, is equal to
$\frac{2\delta(\delta+1)(\delta-1)}{3}.$ Hence,
$c_{4;4}^{(2)}(RLT_{n},0)=c_{4}(RLT_{n},0)\times c_{4}(RLT_{n},0)
-2A+2B-c_{4}(RLT_{n},v)=$
$$\frac{2(\delta-1)\delta(\delta+1)}{3}\times
\Bigl(\frac{2(\delta-1)\delta(\delta+1)}{3}-1\Bigr)
-\frac{(\delta-1)\delta(\delta+1)(13\delta^{2}-7\delta-10)}{6}+$$
$$\frac{(\delta-1)\delta(\delta+1)(27\delta-22)}{12}=
\frac{(\delta-2)(\delta-1)\delta(\delta+1)(2\delta-5)(8\delta-3)}{36}.$$

By duality, we have $c_{5;3}^{(2)}(RLT_{n},0)=c_{3;5}^{(2)}(RLT_{n},0).$
Thus, $c_{8}^{(2)}(RLT_{n},0)=$

$c_{3;5}^{(2)}(RLT_{n},0)+c_{4;4}^{(2)}(RLT_{n},0)
+c_{5;3}^{(2)}(RLT_{n},0)=$

$$2\times\frac{(\delta-2)(\delta-1)\delta(\delta+1)(30\delta^{2}-82\delta+21)}{120}
+$$
$$+\frac{(\delta-2)(\delta-1)\delta(\delta+1)(2\delta-5)(8\delta-3)}{36}=$$
$$\frac{(\delta-2)(\delta-1)\delta(\delta+1)(85\delta^{2}-238\delta+69)}{90}$$
or, in terms of $n=2\delta+1,$
$$c_{8}^{(2)}(RLT_{n},0)=\frac{(n+1)(n-1)(n-3)(n-5)(85n^{2}-646n+837)}{5760}.$$

\bigskip

\noindent {\bf Acknowledgements}

\smallskip

\noindent The author thanks A. Prokofiev for his help in the computer search.

\bigskip

{\bf APPENDIX A: ON $\bold{8}$-CYCLES IN SOME REGULAR TOURNAMENTS OF SMALL ORDERS}

\smallskip

{\bf Appendix A.1: The unique maximizer of $\bold{c_{8}(T)}$ in the class of
regular tournaments of order $\bold{13}$}

\smallskip

Our computer processing of B. McKay's file of regular tournaments of order $13$
implies that the maximum number of $8$-cycles in this class is attained only at
the nearly-doubly-regular tournament $SNDR_{13}$ on the vertex-set
$\{1,2,3,4,5,6,7,8,9,$ $10,11,12,13\}$ with the adjacency matrix
$$\pmatrix 0 & 1 & 1 & 1 & 1 & 1 & 1 & 0 & 0 & 0 & 0 & 0 & 0\\
           0 & 0 & 1 & 0 & 1 & 0 & 1 & 1 & 0 & 0 & 0 & 1 & 1\\
           0 & 0 & 0 & 1 & 0 & 1 & 1 & 0 & 0 & 0 & 1 & 1 & 1\\
           0 & 1 & 0 & 0 & 0 & 1 & 1 & 0 & 1 & 1 & 0 & 0 & 1\\
           0 & 0 & 1 & 1 & 0 & 0 & 0 & 1 & 1 & 1 & 1 & 0 & 0\\
           0 & 1 & 0 & 0 & 1 & 0 & 0 & 0 & 1 & 0 & 1 & 1 & 1\\
           0 & 0 & 0 & 0 & 1 & 1 & 0 & 1 & 0 & 1 & 1 & 1 & 0\\
           1 & 0 & 1 & 1 & 0 & 1 & 0 & 0 & 1 & 0 & 0 & 1 & 0\\
           1 & 1 & 1 & 0 & 0 & 0 & 1 & 0 & 0 & 1 & 1 & 0 & 0\\
           1 & 1 & 1 & 0 & 0 & 1 & 0 & 1 & 0 & 0 & 0 & 0 & 1\\
           1 & 1 & 0 & 1 & 0 & 0 & 0 & 1 & 0 & 1 & 0 & 1 & 0\\
           1 & 0 & 0 & 1 & 1 & 0 & 0 & 0 & 1 & 1 & 0 & 0 & 1\\
           1 & 0 & 0 & 0 & 1 & 0 & 1 & 1 & 1 & 0 & 1 & 0 & 0\endpmatrix.$$

According to our computer search, the distribution of $8$-cycles
among the vertices looks as follows

$c_{8}(SNDR_{13},1)=30618;$

$c_{8}(SNDR_{13},2)=30604;$

$c_{8}(SNDR_{13},3)=30610;$

$c_{8}(SNDR_{13},4)=30598;$

$c_{8}(SNDR_{13},5)=30618;$

$c_{8}(SNDR_{13},6)=30608;$

$c_{8}(SNDR_{13},7)=30604;$

$c_{8}(SNDR_{13},8)=30594;$

$c_{8}(SNDR_{13},9)=30594;$

$c_{8}(SNDR_{13},10)=30612;$

$c_{8}(SNDR_{13},11)=30612;$

$c_{8}(SNDR_{13},12)=30598;$

$c_{8}(SNDR_{13},13)=30610.$

Note that at most two vertices of $SNDR_{13}$ lie on the same number of
$8$-cycles. However, for an automorphism $\pi$ of
a tournament, the length of a $\pi$-orbit
cannot be equal to $2.$ Hence, the automorphism group of $SNDR_{13}$ is trivial.
According to Mathematica, the characteristic polynomial of the adjacency matrix
of $SNDR_{13}$ is $4434+9749x+18310x^{2}+20147x^{3}+19749x^{4}+
13408x^{5}+8358x^{6}+3597x^{7}+1482x^{8}+351x^{9}+91x^{10}-x^{13}.$
The discriminant of this polynomial of degree $13$ is equal to
$$157525764385770965120257003012282911852530325>0.$$
Thus, the spectrum of $SNDR_{13}$ is simple.

\bigskip

{\bf Appendix A.2: The minimizers of $\bold{c_{8}(T)}$ in the classes of
regular tournaments of orders $\bold{n=9,11,13}$}

\smallskip

Let $F_{n}$ be the map which sends a tournament $T$ of order $n$ with the adjacency matrix $A$
to a regular tournament $F_{n}(T)$ of order $2n+1$ with the adjacency matrix
$$F_{n}(A)=\pmatrix A & A^{\top} & \vec{1}\\
	    A^{\top}+I & A       & \vec{0}\\
	   \vec{0}^{\top} & \vec{1}^{\top} & 0\endpmatrix,$$
where $A^{\top}$ is the transpose of $A,$ $I$ is the identity matrix of order $n,$
$\vec{1}$ is the all ones vector of order $n,$ and $\vec{0}$ is the zero vector of order $n.$
It is well known that $F_{n}(TT_{n})=RLT_{2n+1},$ where $TT_{n}$ is the
transitive tournament of order $n.$

Our computer processing of B. McKay's files of regular tournaments of orders $9,11,13$
(see http://cs.anu.edu.au/$\sim$ bdm/data/digraphs.html)
implies that the minimum of $c_{8}(T)$ in these classes is attained only at
$\Delta\circ\Delta,$
$F_{5}\bigl(\Delta(\bullet\Rightarrow\bullet,\bullet,\bullet)\Rightarrow \bullet\bigr),$
and $F_{6}\bigl(\Delta(\bullet\Rightarrow\bullet,\bullet\Rightarrow\bullet,\bullet)\Rightarrow \bullet\bigr),$
respectively.\footnote[3]
{Unfortuntely, $\Delta\circ\Delta$ is not isomorphic to
$F_{4}\bigl(\Delta(\bullet,\bullet,\bullet)\Rightarrow \bullet\bigr).$}
Denote these regular tournaments
by $Umin_{9},Umin_{11},$ and $Umin_{13},$ respectively.
In sequel, we also consider the unique maximizers
$Umax_{9},Umax_{11},$ and $Umax_{13}$ of $c_{8}(T)$
in the classes $\Cal{R}_{9},$ $\Cal{R}_{11},$
and $\Cal{R}_{13},$ respectively.
Recall that $Umax_{9}$ is the unique
vertex-transitive nearly-doubly-regular tournament of order $9,$
$Umax_{11}=QR_{11},$ and the maximizer $Umax_{13}$
has been described above. The computer search yields
$$405=c_{8}(Umin_{9})<   441=c_{8}(RLT_{9}) <        c_{8}(Umax_{9})=477,$$
$$6605=c_{8}(Umin_{11})<   6644=c_{8}(RLT_{11}) <        c_{8}(Umax_{11})=7425,$$
$$45475=c_{8}(Umin_{13})<   45903=c_{8}(RLT_{13}) <        c_{8}(Umax_{13})=49735.$$

We see that $c_{8}(RLT_{9})$ is the arithmetical mean of
the corresponding extreme values $c_{8}(Umin_{9})$ and
$c_{8}(Umax_{9}).$
However, for $n=11,13,$ the value of
$c_{8}(RLT_{n})$ is much closer to the minimum of $c_{8}(T)$ than to
its maximum.

Let $n\overline{c}_{3}^{+}(T)=\sum\limits_{i=1}^{n}c_{3}\bigl(N^{+}(i)\bigr)$
and
$tr_{4}(T)$ be the trace of the $4$th power of the adjacency matrix $A$ of $T.$
Then the formula for $c_{4}(T)$ given in the introduction can be written
in the following form:
$$\overline{c}_{3}^{+}(T)=\frac{(n+1)(n-1)(n-3)}{48}-
\frac{tr_{4}(T)}{4n}.$$

Raising the adjacency matrices of the minimisers $Umin_{11}$ and $Umin_{13}$
to the $4$th power yields $tr_{4}(Umin_{11})$ $=784$ and
$tr_{4}(Umin_{13})=1628.$ Hence,
$$\overline{c}_{3}^{+}(Umin_{11})=
\frac{12*10*8}{48}-17\frac{9}{11}=20-17\frac{9}{11}=2\frac{2}{11}<\frac{5}{2}
=\frac{\overline{c}_{3}^{+}(DR_{11})}{2}$$
and
$$\overline{c}_{3}^{+}(Umin_{13})=
\frac{14*12*10}{48}-31\frac{4}{13}=35-31\frac{4}{13}=3\frac{9}{13}<
4=\frac{\overline{c}_{3}^{+}(NDR_{13})}{2}.$$
Thus, the local properties of $Umin_{11}$ and $Umin_{13}$
are closer to those of $RLT_{11}$ and $RLT_{13}$ than to those of a tournament
with regular structure of the out-sets and in-sets of its vertices.

\bigskip

\centerline{\bf References }

\smallskip

[1]  B. Alspach and C. Tabib, A note on the number of $4$-circuits in a tournament,
Ann. Discrete Math. 12 (1982), 13-19.

[2] A. Asti\'e-Vidal, V. Dugat, Autonomous parts and decomposition of
regular tournaments, Discrete Math. 111 (1993), 27-36.

[3] D. Berman, On the number of $5$-cycles in a tournament,
Congress. Num. 16 (1976), 101-108.

[4] U. Colombo, Sui circuiti nei grafi completi,
Boll. Un. Mat. Ital. 19 (1964), 153-170.

[5] A. Grzesik, D. Kr\'al', L.M. Lov\'asz, J. Volec,
Cycles of a given length in tournaments,
J. Combin. Theory Ser. B 158 (2023), 117-145.

[6] M.G. Kendall and B. Babington Smith, On the method of
paired comparisons, Biometrika 33 (1940), 239-251.

[7] N. Komarov, J. Mackey, On the number of $5$-cycles in a tournament,
J. Graph Theory 86 (2017), 341-356.

[8] A. Kotzig, Sur le nombre des $4$-cycles dans un tournoi,
Mat. Casopis Sloven. Akad. Vied. 18 (1968), 247-254.

[9] J. Plesnik, On homogeneous tournaments, Acta Fac. Rerum Natur.
Univ. Comenian Math. Publ. 21 (1968), 25-34.

[10] K.B. Reid and E. Brown, Doubly-regular tournaments are equivalent to skew-Hada- mard
matrices, J. Combin. Theory Ser. A 12 (1972), 332-338.

[11] P. Rowlison, On $4$-cycles and $5$-cycles in regular tournaments,
Bull. London Math. Soc. 18 (1986), 135-139.

[12] S.V. Savchenko, On $5$-cycles and $6$-cycles in regular $n$-tournaments,
J. Graph Theory 83 (2016), 44-77.

[13] S.V. Savchenko, On the number of $7$-cycles in a regular
$n$-tournament, Discrete Mathematics 340 (2017), 264-285.

[14] S.V. Savchenko, Bernoulli numbers and $m$-cycles
in regular $n$-tournaments, Proc. Amer. Math. Soc.
(will be submitted for publication).

[15] S.V. Savchenko, On $5$-cycles and strong $5$-subtournaments
in a tournament of odd order $n,$ J. Graph Theory
(submitted for publication).

[16] C. Tabib, The number of $4$-cycles in regular tournaments, Utilitas Math. 22 (1982), 315-322.

\end{document}